\documentclass{article}%
\usepackage{amssymb}
\usepackage{amsmath}
\usepackage{makeidx}
\usepackage{amsfonts}
\usepackage{graphicx}
\usepackage{cite}%
\begin{document}

\author{Hartmut Wachter\thanks{E-Mail: Hartmut.Wachter@gmx.de}\\An der Schafscheuer 56\\D-91781 Wei\ss enburg, Federal Republic of Germany}
\title{Conservation laws for a $q$-deformed nonrel\-a\-tivis\-tic particle}
\maketitle
\date{}

\begin{abstract}
We derive $q$-versions of Green's theorem from the Leibniz rules of
partial derivatives for the $q$-deformed Euclidean space. Using these results
and the Schr\"{o}dinger equations for a $q$-deformed nonrelativistic particle,
we derive continuity equations for the probability density, the energy
density, and the momentum density of a $q$-deformed nonrelativistic particle.

\end{abstract}

\section{Introduction}

As is well-known, the divergences that plagued QED from the beginning were
eliminated by the concept of renormalization. On the other hand, various
attempts to overcome the problem with divergences in QED by introducing a
fundamental length have not yet been successful \cite{Hagar:2014}.
Heisenberg's idea of a "lattice-world", for example, led to a breakdown of
continuous rotational and translational symmetries
\cite{Heisenberg:1930,Heisenberg:1938}.

However, some considerations within the framework of a future theory of
quantum gravity suggest that space-time reveals a discrete structure at small
distances \cite{Garay:1995}. For example, the attempt to increase the accuracy
of a position measurement more and more should disturb the background metric
more and more \cite{Mead:1966zz}. Such a fundamental uncertainty in space
could require a space-time algebra generated by non-com\-mu\-ting coordinates.
Realistic noncommutative space-time algebras can be obtained, for example, by
$q$-deformation \cite{CarowWatamura:1990nk,Faddeev:1987ih}.

The $q$-de\-formed Euclidean space has symmetries we can interpret as
$q$-ana\-logs of continuous rotational and translational symmetries. Due to
this property, there are conservation laws for energy and momentum of a
nonrelativistic particle living in $q$-de\-formed Euclidean space. In what
follows, we are briefly describing how we derive these conservation laws in
this article.

First of all, we summarize some results and characteristic features of our
approach [cf. Chap.~\ref{PrelKap}]. The commutation relations for the
coordinates of $q$-de\-formed Euclidean space satisfy the so-called
Poincar\'{e}-Birkhoff-Witt property. Due to this fact, we can associate the
noncommutative algebra of $q$-de\-formed Euclidean space with a commutative
coordinate algebra by using the star-prod\-uct formalism \cite{Moyal:1949sk}.
The star-prod\-uct formalism enables us to construct a $q$-de\-formed version
of mathematical analysis \cite{Carnovale:1999,Wachter:2007A}. As was shown in
Ref.~\cite{Wachter:2020A}, the time evolution operator of a quantum system in
$q$-de\-formed Euclidean space is of the same form as in the undeformed case.
These findings enable us to write down Schr\"{o}\-dinger equations for a
$q$-de\-formed nonrelativistic particle (cf. Ref.~\cite{Wachter:2020B} and
Chap.~\ref{KapHerSchGle}).

In Chap.~\ref{KapGre} of this article, we derive $q$-versions of Green's
theorem. We achieve this by applying the Leibniz rules for
$q$-de\-formed partial derivatives. From the $q$-versions of Green's theorem and the Schr\"{o}\-dinger equations for a $q$-de\-formed particle, we
derive continuity equations for the probability density, the energy density,
and the momentum density of a $q$-de\-formed, nonrelativistic particle. Our
reasonings also include the case where a $q$-de\-formed particle interacts
with an electromagnetic field. In this respect, we will show
in\ Chap.~\ref{EicTraKap} that our $q$-de\-formed continuity equations are
invariant under gauge transformations. Finally, we calculate Heisenberg's
equation for observables of a $q$-de\-formed nonrelativistic particle
interacting with an electromagnetic field. In doing so, we obtain evolution
equations that are in agreement with our $q$-de\-formed continuity equations.

\section{Preliminaries\label{PrelKap}}

\subsection{Star-products\label{KapQuaZeiEle}}

The three-di\-men\-sion\-al $q$-de\-formed Euclidean space $\mathbb{R}_{q}%
^{3}$ has the generators $X^{+}$, $X^{3}$, and $X^{-}$, subject to the
following commutation relations \cite{Lorek:1997eh}:%
\begin{align}
X^{3}X^{+} &  =q^{2}X^{+}X^{3},\nonumber\\
X^{3}X^{-} &  =q^{-2}X^{-}X^{3},\nonumber\\
X^{-}X^{+} &  =X^{+}X^{-}+(q-q^{-1})\hspace{0.01in}X^{3}X^{3}%
.\label{RelQuaEukDre}%
\end{align}
We can extend the algebra of $\mathbb{R}_{q}^{3}$ by a time element $X^{0}$,
which commutes with the generators $X^{+}$, $X^{3}$, and $X^{-}$
\cite{Wachter:2020A}:%
\begin{equation}
X^{0}X^{A}=X^{A}X^{0},\text{\qquad}A\in\{+,3,-\}.\label{ZusRelExtDreEukQUa}%
\end{equation}
In the following, we refer to the algebra spanned by the generators $X^{i}$
with $i\in\{0,+,3,-\}$ as $\mathbb{R}_{q}^{3,t}$.

There is a $q$-ana\-log of the three-di\-men\-sion\-al Euclidean metric
$g^{AB}$ with its inverse $g_{AB}$ \cite{Lorek:1997eh} (rows and columns are
arranged in the order $+,3,-$):%
\begin{equation}
g_{AB}=g^{AB}=\left(
\begin{array}
[c]{ccc}%
0 & 0 & -\hspace{0.01in}q\\
0 & 1 & 0\\
-\hspace{0.01in}q^{-1} & 0 & 0
\end{array}
\right)  .\label{DreDimMet}%
\end{equation}
We can use the $q$-de\-formed metric to raise and lower indices:%
\begin{equation}
X_{A}=g_{AB}\hspace{0.01in}X^{B},\qquad X^{A}=g^{AB}X_{B}.\label{HebSenInd}%
\end{equation}

The algebra $\mathbb{R}_{q}^{3,t}$ has a semilinear, involutive, and
anti-multiplicative mapping, which we call \textit{quantum space conjugation}.
If we indicate conjugate elements of a quantum space by a bar,\footnote{A bar
over a complex number indicates complex conjugation.} we can write the
properties of quantum space conjugation as follows ($\alpha,\beta\in
\mathbb{C}$ and $u,v\in\mathbb{R}_{q}^{3,t}$):%
\begin{equation}
\overline{\alpha\,u+\beta\,v}=\overline{\alpha}\,\overline{u}+\overline{\beta
}\,\overline{v},\quad\overline{\overline{u}}=u,\quad\overline{u\,v}%
=\overline{v}\,\overline{u}.
\end{equation}
The conjugation for $\mathbb{R}_{q}^{3,t}$ is compatible with the commutation
relations in Eq.~(\ref{RelQuaEukDre}) and Eq.~(\ref{ZusRelExtDreEukQUa}) if
the following applies \cite{Wachter:2020A}:%
\begin{equation}
\overline{X^{A}}=X_{A}=g_{AB}\hspace{0.01in}X^{B},\qquad\overline{X^{0}}%
=X_{0}. \label{ConSpaKoo}%
\end{equation}

We can only prove a physical theory if it predicts measurement results. The
problem, however, is: How can we associate the elements of the noncommutative
space $\mathbb{R}_{q}^{3,t}$ with real numbers? One solution to this problem
is to introduce a vector space isomorphism between the noncommutative algebra
$\mathbb{R}_{q}^{3,t}$ and a corresponding commutative coordinate algebra
$\mathbb{C}[\hspace{0.01in}x^{+},x^{3},x^{-},t\hspace{0.01in}]$.

We recall that the nor\-mal-or\-dered monomials in the generators $X^{i}$ form
a basis of the algebra $\mathbb{R}_{q}^{3,t}$, i.~e. we can write each element
$F\in$ $\mathbb{R}_{q}^{3,t}$ uniquely as a finite or infinite linear
combination of monomials with a given normal ordering
(\textit{Poincar\'{e}-Birkhoff-Witt property}):%
\begin{equation}
F=\sum\limits_{n_{+},\ldots,\hspace{0.01in}n_{0}}a_{\hspace{0.01in}n_{+}%
\ldots\hspace{0.01in}n_{0}}\,(X^{+})^{n_{+}}(X^{3})^{n_{3}}(X^{-})^{n_{-}%
}(X^{0})^{n_{0}},\quad\quad a_{\hspace{0.01in}n_{+}\ldots\hspace{0.01in}n_{0}%
}\in\mathbb{C}.
\end{equation}
Since the monomials $(x^{+})^{n_{+}}(x^{3})^{n_{3}}(x^{-})^{n_{-}}%
\hspace{0.01in}t^{\hspace{0.01in}n_{0}}$ with $n_{+},\ldots,n_{0}\in
\mathbb{N}_{0}$ form a basis of the commutative algebra $\mathbb{C}%
[\hspace{0.01in}x^{+},x^{3},x^{-},t\hspace{0.01in}]$, we can define a vector
space isomorphism%
\begin{equation}
\mathcal{W}:\mathbb{C}[\hspace{0.01in}x^{+},x^{3},x^{-},t\hspace
{0.01in}]\rightarrow\mathbb{R}_{q}^{3,t} \label{VecRauIsoInv}%
\end{equation}
with%
\begin{equation}
\mathcal{W}\left(  (x^{+})^{n_{+}}(x^{3})^{n_{3}}(x^{-})^{n_{-}}%
\hspace{0.01in}t^{\hspace{0.01in}n_{0}}\right)  =(X^{+})^{n_{+}}(X^{3}%
)^{n_{3}}(X^{-})^{n_{-}}(X^{0})^{n_{0}}. \label{StePro0}%
\end{equation}
In general, we have%
\begin{equation}
\mathbb{C}[\hspace{0.01in}x^{+},x^{3},x^{-},t\hspace{0.01in}]\ni f\mapsto
F\in\mathbb{R}_{q}^{3,t},
\end{equation}
where%
\begin{align}
f  &  =\sum\limits_{n_{+},\ldots,\hspace{0.01in}n_{0}}a_{\hspace{0.01in}%
n_{+}\ldots\hspace{0.01in}n_{0}}\,(x^{+})^{n_{+}}(x^{3})^{n_{3}}(x^{-}%
)^{n_{-}}\hspace{0.01in}t^{\hspace{0.01in}n_{0}},\nonumber\\
F  &  =\sum\limits_{n_{+},\ldots,\hspace{0.01in}n_{0}}a_{\hspace{0.01in}%
n_{+}\ldots\hspace{0.01in}n_{0}}\,(X^{+})^{n_{+}}(X^{3})^{n_{3}}(X^{-}%
)^{n_{-}}(X^{0})^{n_{0}}. \label{AusFfNorOrd}%
\end{align}
The vector space isomorphism $\mathcal{W}$ is nothing else but the
\textit{Moy\-al-Weyl mapping}, which gives an operator $F$ to a complex valued
function $f$
\cite{Bayen:1977ha,1997q.alg.....9040K,Madore:2000en,Moyal:1949sk}.

We can extend this vector space isomorphism to an algebra isomorphism if we
introduce a new product on the commutative coordinate algebra. This so-called
\textit{star-pro\-duct }symbolized by $\circledast$ satisfies the following
homomorphism condition:%
\begin{equation}
\mathcal{W}\left(  f\circledast g\right)  =\mathcal{W}\left(  f\right)
\cdot\mathcal{W}\left(  \hspace{0.01in}g\right)  . \label{HomBedWeyAbb}%
\end{equation}
Since the Moy\-al-Weyl mapping is invertible, we can write the star-prod\-uct
as follows:%
\begin{equation}
f\circledast g=\mathcal{W}^{\hspace{0.01in}-1}\big (\,\mathcal{W}\left(
f\right)  \cdot\mathcal{W}\left(  \hspace{0.01in}g\right)  \big ).
\label{ForStePro}%
\end{equation}

To get explicit formulas for calculating star-prod\-ucts, we first have to
write a noncommutative product of two nor\-mal-or\-dered monomials as a linear
combination of nor\-mal-or\-dered monomials again (see
Ref.~\cite{Wachter:2002A} for details):%
\begin{equation}
(X^{+})^{n_{+}}\ldots\hspace{0.01in}(X^{0})^{n_{0}}\cdot(X^{+})^{m_{+}}%
\ldots\hspace{0.01in}(X^{0})^{m_{0}}=\sum_{\underline{k}\hspace{0.01in}%
=\hspace{0.01in}0}B_{\underline{k}}^{\hspace{0.01in}\underline{n}%
,\underline{m}}\,(X^{+})^{k_{+}}\ldots\hspace{0.01in}(X^{0})^{k_{0}}.
\label{EntProMon}%
\end{equation}
We achieve this by using the commutation relations for the noncommutative
coordinates [cf. Eq.~(\ref{RelQuaEukDre})]. From the concrete form of the
expansion in\ Eq.~(\ref{EntProMon}), we can finally read off a formula to
calculate the star-prod\-uct of two power series in commutative space-time
coordinates ($\lambda=q-q^{-1}$):\footnote{The argument $\mathbf{x}$ indicates
a dependence on the spatial coordinates $x^{+}$, $x^{3}$, and $x^{-}$.}%
\begin{gather}
f(\mathbf{x},t)\circledast g(\mathbf{x},t)=\nonumber\\
\sum_{k\hspace{0.01in}=\hspace{0.01in}0}^{\infty}\lambda^{k}\hspace
{0.01in}\frac{(x^{3})^{2k}}{[[k]]_{q^{4}}!}\,q^{2(\hat{n}_{3}\hspace
{0.01in}\hat{n}_{+}^{\prime}+\,\hat{n}_{-}\hat{n}_{3}^{\prime})}%
D_{q^{4},\hspace{0.01in}x^{-}}^{k}f(\mathbf{x},t)\,D_{q^{4},\hspace
{0.01in}x^{\prime+}}^{k}g(\mathbf{x}^{\prime},t)\big|_{x^{\prime}%
\rightarrow\hspace{0.01in}x}. \label{StaProForExp}%
\end{gather}
The expression\ above depends on the operators%
\begin{equation}
\hat{n}_{A}=x^{A}\frac{\partial}{\partial x^{A}}%
\end{equation}
and the so-called Jackson derivatives \cite{Jackson:1910yd}:%
\begin{equation}
D_{q^{k},\hspace{0.01in}x}\,f=\frac{f(q^{k}x)-f(x)}{q^{k}x-x}.
\end{equation}
Moreover, the $q$-numbers are given by%
\begin{equation}
\lbrack\lbrack a]]_{q}=\frac{1-q^{a}}{1-q},
\end{equation}
and the $q$-factorials are defined in complete analogy to the undeformed case:%
\begin{equation}
\lbrack\lbrack n]]_{q}!=[[1]]_{q}\hspace{0.01in}[[2]]_{q}\ldots\lbrack
\lbrack\hspace{0.01in}n-1]]_{q}\hspace{0.01in}[[n]]_{q},\qquad\lbrack
\lbrack0]]_{q}!=1.
\end{equation}

The algebra isomorphism $\mathcal{W}^{-1}$ also enables us to carry over the
conjugation for the quantum space algebra $\mathbb{R}_{q}^{3,t}$ to the
commutative coordinate algebra $\mathbb{C}[\hspace{0.01in}x^{+},x^{3}%
,x^{-},t\hspace{0.01in}]$. In other words, the mapping $\mathcal{W}%
^{\hspace{0.01in}-1}$ is a $\ast$-al\-ge\-bra homomorphism:%
\begin{equation}
\mathcal{W}(\hspace{0.01in}\overline{f}\hspace{0.01in})=\overline
{\mathcal{W}(f)}\qquad\Leftrightarrow\text{\qquad}\overline{f}=\mathcal{W}%
^{-1}\big (\hspace{0.01in}\overline{\mathcal{W}(f)}\hspace{0.01in}\big ).
\label{ConAlgIso}%
\end{equation}
This relationship implies the following property for the star-pro\-duct:%
\begin{equation}
\overline{f\circledast g}=\overline{g}\circledast\overline{f}.
\label{KonEigSteProFkt}%
\end{equation}

With $\bar{f}$, we designate the power series obtained from $f$ by quantum
space conjugation. It follows from Eq.~(\ref{ConSpaKoo})\ and
Eq.~(\ref{ConAlgIso}) that $\bar{f}$ takes the following form (if $\bar
{a}_{n_{+},n_{3},n_{-},n_{0}}$ stands for the complex conjugate of
$a_{n_{+},n_{3},n_{-},n_{0}}$) \cite{Wachter:2007A,Wachter:2020A}:%
\begin{align}
\overline{f(\mathbf{x},t)}  &  =\sum\nolimits_{\underline{n}}\bar{a}%
_{n_{+},n_{3},n_{-},n_{0}}\,(-\hspace{0.01in}q\hspace{0.01in}x^{-})^{n_{+}%
}(\hspace{0.01in}x^{3})^{n_{3}}(-\hspace{0.01in}q^{-1}x^{+})^{n_{-}}%
\hspace{0.01in}t^{\hspace{0.01in}n_{0}}\nonumber\\
&  =\sum\nolimits_{\underline{n}}(-\hspace{0.01in}q)^{n_{-}-\hspace
{0.02in}n_{+}}\hspace{0.01in}\bar{a}_{n_{-},n_{3},n_{+},n_{0}}\,(\hspace
{0.01in}x^{+})^{n_{+}}(\hspace{0.01in}x^{3})^{n_{3}}(\hspace{0.01in}%
x^{-})^{n_{-}}\hspace{0.01in}t^{\hspace{0.01in}n_{0}}\nonumber\\
&  =\bar{f}(\mathbf{x},t). \label{KonPotReiKom}%
\end{align}

\subsection{Partial derivatives and integrals\label{KapParDer}}

There are partial derivatives for $q$-de\-formed space-time coordinates
\cite{CarowWatamura:1990zp,Wess:1990vh}. These partial derivatives again form
a quantum space with the same algebraic structure as that of the
$q$-de\-formed space-time coordinates. Thus, the $q$-de\-formed partial
derivatives $\partial_{i}$ satisfy the same commutation relations as the
covariant coordinate generators $X_{i}$:%
\begin{gather}
\partial_{0}\hspace{0.01in}\partial_{+}=\hspace{0.01in}\partial_{+}%
\hspace{0.01in}\partial_{0},\quad\partial_{0}\hspace{0.01in}\partial
_{-}=\hspace{0.01in}\partial_{-}\hspace{0.01in}\partial_{0},\quad\partial
_{0}\hspace{0.01in}\partial_{3}=\partial_{3}\hspace{0.01in}\partial
_{0},\nonumber\\
\partial_{+}\hspace{0.01in}\partial_{3}=q^{2}\partial_{3}\hspace
{0.01in}\partial_{+},\quad\partial_{3}\hspace{0.01in}\partial_{-}%
=\hspace{0.01in}q^{2}\partial_{-}\hspace{0.01in}\partial_{3},\nonumber\\
\partial_{+}\hspace{0.01in}\partial_{-}-\partial_{-}\hspace{0.01in}%
\partial_{+}=\hspace{0.01in}\lambda\hspace{0.01in}\partial_{3}\hspace
{0.01in}\partial_{3}.
\end{gather}
The commutation relations above are invariant under conjugation if the
derivatives show the following conjugation properties:\footnote{The indices of
partial derivatives are raised and lowered in the same way as those of
coordinates [see Eq.~(\ref{HebSenInd}) in Chap.~\ref{KapQuaZeiEle}].}%
\begin{equation}
\overline{\partial_{A}}=-\hspace{0.01in}\partial^{A}=-g^{AB}\partial
_{B},\qquad\overline{\partial_{0}}=-\hspace{0.01in}\partial^{0}=-\hspace
{0.01in}\partial_{0}. \label{KonAbl}%
\end{equation}

There are two ways of commuting $q$-de\-formed partial derivatives with
$q$-de\-formed space-time coordinates. One is given by the following
$q$-de\-formed Leibniz rules
\cite{CarowWatamura:1990zp,Wess:1990vh,Wachter:2020A}:%
\begin{align}
\partial_{B}X^{A} &  =\delta_{B}^{A}+q^{4}\hat{R}{^{AC}}_{BD}\,X^{D}%
\partial_{C},\nonumber\\
\partial_{A}X^{0} &  =X^{0}\hspace{0.01in}\partial_{A},\nonumber\\
\partial_{0}\hspace{0.01in}X^{A} &  =X^{A}\hspace{0.01in}\partial
_{0},\nonumber\\
\partial_{0}\hspace{0.01in}X^{0} &  =1+X^{0}\hspace{0.01in}\partial
_{0}.\label{DifKalExtEukQuaDreUnk}%
\end{align}
Note that $\hat{R}{^{AC}}_{BD}$ denotes the vector representation of the
R-ma\-trix for the three-di\-men\-sion\-al $q$-de\-formed Euclidean space.

By conjugation, we can obtain the Leibniz rules for another differential
calculus from the identities in Eq.~(\ref{DifKalExtEukQuaDreUnk}). Introducing
$\hat{\partial}_{A}=q^{6}\partial_{A}$ and $\hat{\partial}_{0}=\partial_{0}$,
we can write the Leibniz rules of this second differential calculus in the
following form:%
\begin{align}
\hat{\partial}_{B}\hspace{0.01in}X^{A}  &  =\delta_{B}^{A}+q^{-4}(\hat{R}%
^{-1}){^{AC}}_{BD}\,X^{D}\hat{\partial}_{C},\nonumber\\
\hat{\partial}_{A}\hspace{0.01in}X^{0}  &  =X^{0}\hspace{0.01in}\hat{\partial
}_{A},\nonumber\\
\hat{\partial}_{0}\hspace{0.01in}X^{A}  &  =X^{A}\hspace{0.01in}\hat{\partial
}_{0},\nonumber\\
\hat{\partial}_{0}\hspace{0.01in}X^{0}  &  =1+X^{0}\hspace{0.01in}%
\hat{\partial}_{0}. \label{DifKalExtEukQuaDreKon}%
\end{align}

Using the Leibniz rules in Eq.$~$(\ref{DifKalExtEukQuaDreUnk}) or
Eq.$~$(\ref{DifKalExtEukQuaDreKon}), we can calculate how partial derivatives
act on nor\-mal-or\-dered monomials of noncommutative coordinates. We can
carry over these actions to commutative coordinate monomials with the help of
the Mo\-yal-Weyl mapping:%
\begin{equation}
\partial^{i}\triangleright(x^{+})^{n_{+}}(x^{3})^{n_{3}}(x^{-})^{n_{-}}%
\hspace{0.01in}t^{\hspace{0.01in}n_{0}}=\mathcal{W}^{\hspace{0.01in}%
-1}\big (\partial^{i}\triangleright(X^{+})^{n_{+}}(X^{3})^{n_{3}}%
(X^{-})^{n_{-}}(X^{0})^{n_{0}}\big ).
\end{equation}
Since the Mo\-yal-Weyl mapping is linear, we can apply the action above to
space-time functions that can be written as a power series:%
\begin{equation}
\partial^{i}\triangleright f(\mathbf{x},t)=\mathcal{W}^{\hspace{0.01in}%
-1}\big (\partial^{i}\triangleright\mathcal{W}(f(\mathbf{x},t))\big ).
\end{equation}

If we use the ordering given in Eq.~(\ref{StePro0}) of the previous chapter,
the Leibniz rules in Eq.~(\ref{DifKalExtEukQuaDreUnk})\ will lead to the
following operator representations \cite{Bauer:2003}:%
\begin{align}
\partial_{+}\triangleright f(\mathbf{x},t)  &  =D_{q^{4},\hspace{0.01in}x^{+}%
}f(\mathbf{x},t),\nonumber\\
\partial_{3}\triangleright f(\mathbf{x},t)  &  =D_{q^{2},\hspace{0.01in}x^{3}%
}f(q^{2}x^{+},x^{3},x^{-},t),\nonumber\\
\partial_{-}\triangleright f(\mathbf{x},t)  &  =D_{q^{4},\hspace{0.01in}x^{-}%
}f(x^{+},q^{2}x^{3},x^{-},t)+\lambda\hspace{0.01in}x^{+}D_{q^{2}%
,\hspace{0.01in}x^{3}}^{2}f(\mathbf{x},t). \label{UnkOpeDarAbl}%
\end{align}
The derivative $\partial_{0}$, however, is represented on the commutative
space-time algebra by an ordinary partial derivative:%
\begin{equation}
\partial_{0}\triangleright\hspace{-0.01in}f(\mathbf{x},t)=\partial
_{t}\triangleright\hspace{-0.01in}f(\mathbf{x},t)=\frac{\partial
f(\mathbf{x},t)}{\partial t}. \label{OpeDarZeiAblExtQuaEuk}%
\end{equation}

Using the Leibniz rules in Eq.$~$(\ref{DifKalExtEukQuaDreKon}), we get
operator representations for the partial derivatives $\hat{\partial}_{i}$. The
Leibniz rules in Eq.$~$(\ref{DifKalExtEukQuaDreUnk}) and Eq.$~$%
(\ref{DifKalExtEukQuaDreKon}) are transformed into each other by the following
substitutions:%
\begin{gather}
q\rightarrow q^{-1},\quad X^{-}\rightarrow X^{+},\quad X^{+}\rightarrow
X^{-},\nonumber\\
\partial^{\hspace{0.01in}+}\rightarrow\hat{\partial}^{\hspace{0.01in}-}%
,\quad\partial^{\hspace{0.01in}-}\rightarrow\hat{\partial}^{\hspace{0.01in}%
+},\quad\partial^{\hspace{0.01in}3}\rightarrow\hat{\partial}^{\hspace
{0.01in}3},\quad\partial^{\hspace{0.01in}0}\rightarrow\hat{\partial}%
^{\hspace{0.01in}0}. \label{UebRegGedUngAblDreQua}%
\end{gather}
For this reason, we obtain the operator representations of the partial
derivatives $\hat{\partial}_{A}$ from those of the partial derivatives
$\partial_{A}$ [cf. Eq.~(\ref{UnkOpeDarAbl})] if we replace $q$ by $q^{-1}$
and exchange the indices $+$ and $-$:%
\begin{align}
\hat{\partial}_{-}\,\bar{\triangleright}\,f(\mathbf{x},t)  &  =D_{q^{-4}%
,\hspace{0.01in}x^{-}}f(\mathbf{x},t),\nonumber\\
\hat{\partial}_{3}\,\bar{\triangleright}\,f(\mathbf{x},t)  &  =D_{q^{-2}%
,\hspace{0.01in}x^{3}}f(q^{-2}x^{-},x^{3},x^{+},t),\nonumber\\
\hat{\partial}_{+}\,\bar{\triangleright}\,f(\mathbf{x},t)  &  =D_{q^{-4}%
,\hspace{0.01in}x^{+}}f(x^{-},q^{-2}x^{3},x^{+},t)-\lambda\hspace{0.01in}%
x^{-}D_{q^{-2},\hspace{0.01in}x^{3}}^{2}f(\mathbf{x},t). \label{KonOpeDarAbl}%
\end{align}
Once again, $\hat{\partial}_{0}$ is represented on the commutative space-time
algebra by an ordinary partial derivative:%
\begin{equation}
\hat{\partial}_{0}\,\bar{\triangleright}\,f(\mathbf{x},t)=\partial
_{t}\triangleright\hspace{-0.01in}f(\mathbf{x},t)=\frac{\partial
f(\mathbf{x},t)}{\partial t}. \label{OpeDarZeiAblExtQuaEukKon}%
\end{equation}
Due to the substitutions given in\ Eq.~(\ref{UebRegGedUngAblDreQua}), the
actions in Eqs.~(\ref{KonOpeDarAbl}) and (\ref{OpeDarZeiAblExtQuaEukKon})
refer to nor\-mal-or\-dered monomials different from those in
Eq.~(\ref{StePro0}) of the previous chapter:%
\begin{equation}
\widetilde{\mathcal{W}}\left(  t^{\hspace{0.01in}n_{0}}(x^{+})^{n_{+}}%
(x^{3})^{n_{3}}(x^{-})^{n_{-}}\right)  =(X^{0})^{n_{0}}(X^{-})^{n_{-}}%
(X^{3})^{n_{3}}(X^{+})^{n_{+}}. \label{UmNor}%
\end{equation}

We should not forget that we can also commute $q$-de\-formed partial
derivatives from the \textit{right} side of a nor\-mal-or\-dered monomial to
the left side by using the Leibniz rules. This way, we get\ so-called
\textit{right} re\-pre\-sen\-ta\-tions of partial derivatives, for which we
write $f\,\bar{\triangleleft}\,\partial^{i}$ or $f\triangleleft\hat{\partial
}^{i}$. Note that the operation of conjugation transforms left actions of
partial derivatives into right actions and vice versa \cite{Bauer:2003}:%
\begin{align}
\overline{\partial^{i}\triangleright f} &  =-\bar{f}\,\bar{\triangleleft
}\,\partial_{i}, & \overline{f\,\bar{\triangleleft}\,\partial^{i}} &
=-\hspace{0.01in}\partial_{i}\triangleright\bar{f},\nonumber\\
\overline{\hat{\partial}^{i}\,\bar{\triangleright}\,f} &  =-\bar
{f}\triangleleft\hat{\partial}_{i}, & \overline{f\triangleleft\hat{\partial
}^{i}} &  =-\hspace{0.01in}\hat{\partial}_{i}\,\bar{\triangleright}\,\bar
{f}.\label{RegConAbl}%
\end{align}

In general, the operator representations in Eqs.~(\ref{UnkOpeDarAbl}) and
(\ref{KonOpeDarAbl}) consist of two terms, which we call $\partial
_{\operatorname*{cla}}^{A}$ and $\partial_{\operatorname*{cor}}^{A}$:%
\begin{equation}
\partial^{A}\triangleright F=\left(  \partial_{\operatorname*{cla}}%
^{A}+\partial_{\operatorname*{cor}}^{A}\right)  \triangleright F.
\end{equation}
In the undeformed limit $q\rightarrow1$, $\partial_{\operatorname*{cla}}^{A}$
becomes an ordinary partial derivative, and $\partial_{\operatorname*{cor}%
}^{A}$ disappears. We get a solution to the difference equation $\partial
^{A}\triangleright F=f$ with given $f$ by using the following formula
\cite{Wachter:2004A}:%
\begin{align}
F  &  =(\partial^{A})^{-1}\triangleright f=\left(  \partial
_{\operatorname*{cla}}^{A}+\partial_{\operatorname*{cor}}^{A}\right)
^{-1}\triangleright f\nonumber\\
&  =\sum_{k\hspace{0.01in}=\hspace{0.01in}0}^{\infty}\left[  -(\partial
_{\operatorname*{cla}}^{A})^{-1}\partial_{\operatorname*{cor}}^{A}\right]
^{k}(\partial_{\operatorname*{cla}}^{A})^{-1}\triangleright f.
\end{align}
Applying the above formula to the operator representations in
Eq.~(\ref{UnkOpeDarAbl}), we get%
\begin{align}
(\partial_{+})^{-1}\triangleright f(\mathbf{x},t)  &  =D_{q^{4},\hspace
{0.01in}x^{+}}^{-1}f(\mathbf{x},t),\nonumber\\
(\partial_{3})^{-1}\triangleright f(\mathbf{x},t)  &  =D_{q^{2},\hspace
{0.01in}x^{3}}^{-1}f(q^{-2}x^{+},x^{3},x^{-},t), \label{InvParAbl1}%
\end{align}
and%
\begin{gather}
(\partial_{-})^{-1}\triangleright f(\mathbf{x},t)=\nonumber\\
=\sum_{k\hspace{0.01in}=\hspace{0.01in}0}^{\infty}q^{2k\left(  k\hspace
{0.01in}+1\right)  }\left(  -\lambda\,x^{+}D_{q^{4},\hspace{0.01in}x^{-}}%
^{-1}D_{q^{2},\hspace{0.01in}x^{3}}^{2}\right)  ^{k}D_{q^{4},\hspace
{0.01in}x^{-}}^{-1}f(x^{+},q^{-2\left(  k\hspace{0.01in}+1\right)  }%
x^{3},x^{-},t). \label{InvParAbl2}%
\end{gather}
Note that $D_{q,\hspace{0.01in}x}^{-1}$ stands for a Jackson integral with $x$
being the variable of integration \cite{Jackson:1908}. The explicit form of
this Jackson integral depends on its limits of integration and the value for
the deformation parameter $q$. If $x>0$ and $q>1$, for example, the following
applies:%
\begin{equation}
\int_{0}^{\hspace{0.01in}x}\text{d}_{q}z\hspace{0.01in}f(z)=(q-1)\hspace
{0.01in}x\sum_{j=1}^{\infty}q^{-j}f(q^{-j}x).
\end{equation}
Finally, the integral for the time coordinate is an ordinary integral since
$\partial_{0}$ acts on the commutative space-time algebra like an ordinary
partial derivative [cf. Eq.~(\ref{OpeDarZeiAblExtQuaEuk})]:%
\begin{equation}
(\partial_{0})^{-1}\triangleright f(\mathbf{x},t)\hspace{0.01in}=\int
\text{d}t\,f(\mathbf{x},t).
\end{equation}

The above considerations also apply to the partial derivatives with a hat.
However, we can obtain the representations of $\hat{\partial}_{i}$ from those
of the derivatives $\partial_{i}$ if we replace $q$ with $q^{-1}$ and exchange
the indices $+$ and $-$. Applying these substitutions to the expressions in
Eqs.~(\ref{InvParAbl1}) and (\ref{InvParAbl2}), we immediately get the
corresponding results for the partial derivatives $\hat{\partial}_{i}$.

By successively applying the integral operators given in
Eqs.~(\ref{InvParAbl1}) and (\ref{InvParAbl2}), we can explain an integration
over all space \cite{Wachter:2004A,Wachter:2007A}:%
\begin{equation}
\int_{-\infty}^{+\infty}\text{d}_{q}^{3}\hspace{0.01in}x\,f(x^{+},x^{3}%
,x^{-})=(\partial_{-})^{-1}\big |_{-\infty}^{+\infty}\,(\partial_{3}%
)^{-1}\big |_{-\infty}^{+\infty}\,(\partial_{+})^{-1}\big |_{-\infty}%
^{+\infty}\triangleright f.
\end{equation}
On the right-hand side of the above relation, the different integral operators
can be simplified to Jackson integrals \cite{Wachter:2004A,Jambor:2004ph}:%
\begin{equation}
\int_{-\infty}^{+\infty}\text{d}_{q}^{3}\hspace{0.01in}x\,f(\mathbf{x}%
)=D_{q^{2},\hspace{0.01in}x^{-}}^{-1}\big |_{-\infty}^{+\infty}\,D_{q,x^{3}%
}^{-1}\big |_{-\infty}^{+\infty}\,D_{q^{2},\hspace{0.01in}x^{+}}%
^{-1}\big |_{-\infty}^{+\infty}\,f(\mathbf{x}).
\end{equation}
Note that the Jackson integrals in the formula above refer to a smaller
$q$-lattice. Using such a smaller $q$-lattice ensures that our integral over
all space is a scalar with trivial braiding properties \cite{Kempf:1994yd}.

The $q$-integral over all space shows some significant features
\cite{Wachter:2007A}. In this respect, $q$-de\-formed versions of
\textit{Stokes' theorem} apply:%
\begin{align}
\int_{-\infty}^{+\infty}\text{d}_{q}^{3}\hspace{0.01in}x\,\partial
^{A}\triangleright f  &  =\int_{-\infty}^{+\infty}\text{d}_{q}^{3}%
\hspace{0.01in}x\,f\,\bar{\triangleleft}\,\partial^{A}=0,\nonumber\\
\int_{-\infty}^{+\infty}\text{d}_{q}^{3}\hspace{0.01in}x\,\hat{\partial}%
^{A}\,\bar{\triangleright}\,f  &  =\int_{-\infty}^{+\infty}\text{d}_{q}%
^{3}\hspace{0.01in}x\,f\triangleleft\hat{\partial}^{A}=0. \label{StoThe}%
\end{align}
The $q$-de\-formed Stokes' theorem also implies rules for integration by
parts:%
\begin{align}
\int_{-\infty}^{+\infty}\text{d}_{q}^{3}\hspace{0.01in}x\,f\circledast
(\partial^{A}\triangleright g)  &  =\int_{-\infty}^{+\infty}\text{d}_{q}%
^{3}\hspace{0.01in}x\,(f\triangleleft\partial^{A})\circledast g,\nonumber\\
\int_{-\infty}^{+\infty}\text{d}_{q}^{3}\hspace{0.01in}x\,f\circledast
(\hat{\partial}^{A}\,\bar{\triangleright}\,g)  &  =\int_{-\infty}^{+\infty
}\text{d}_{q}^{3}\hspace{0.01in}x\,(f\,\bar{\triangleleft}\,\hat{\partial}%
^{A})\circledast g. \label{PatIntUneRaumInt}%
\end{align}
Finally, we mention that the $q$-integral over all space behaves as follows
under quantum space conjugation:%
\begin{equation}
\overline{\int_{-\infty}^{+\infty}\text{d}_{q}^{3}\hspace{0.01in}x\,f}%
=\int_{-\infty}^{+\infty}\text{d}_{q}^{3}\hspace{0.01in}x\,\bar{f}.
\label{KonEigVolInt}%
\end{equation}

\subsection{Hopf structures and L-matrices\label{KapHofStr}}

The three-di\-men\-sio\-nal $q$-de\-formed Euclidean space $\mathbb{R}_{q}%
^{3}$ is a three-di\-men\-sio\-nal representation of the Drin\-feld-Jim\-bo
al\-ge\-bra $\mathcal{U}_{q}(\operatorname*{su}_{2})$. The latter is a
deformation of the universal enveloping algebra of the Lie algebra
$\operatorname*{su}_{2}$ \cite{Kulish:1983md}. Accordingly, the algebra
$\mathcal{U}_{q}(\operatorname*{su}_{2})$ has three generators $T^{+}$,
$T^{-}$, and $T^{3}$ which satisfy the following relations \cite{Lorek:1993tq}%
:%
\begin{align}
q^{-1}\hspace{0.01in}T^{+}T^{-}-q\,T^{-}T^{+}  &  =T^{3},\nonumber\\
q^{\hspace{0.01in}2}\hspace{0.01in}T^{3}T^{+}-q^{-2}\hspace{0.01in}T^{+}T^{3}
&  =(q+q^{-1})\hspace{0.01in}T^{+},\nonumber\\
q^{\hspace{0.01in}2}\hspace{0.01in}T^{-}T^{3}-q^{-2}\hspace{0.01in}T^{3}T^{-}
&  =(q+q^{-1})\hspace{0.01in}T^{-}. \label{VerRel1UqSU2}%
\end{align}

The algebra of the $q$-de\-formed partial derivatives $\partial^{A}$,
$A\in\{+,3,-\}$, together with $\mathcal{U}_{q}(\operatorname*{su}_{2})$ form
the cross-prod\-uct algebra $\mathbb{R}_{q}^{3}\rtimes\hspace{0.01in}%
\mathcal{U}_{q}(\operatorname*{su}\nolimits_{2})$
\cite{majid-1993-34,Weixler:1993ph}. We know that the algebra $\mathbb{R}%
_{q}^{3}\rtimes\mathcal{U}_{q}(\operatorname*{su}\nolimits_{2})$ is a Hopf
algebra \cite{Klimyk:1997eb}. Accordingly, the $q$-de\-formed partial
derivatives as elements of $\mathbb{R}_{q}^{3}\rtimes\hspace{0.01in}%
\mathcal{U}_{q}(\operatorname*{su}\nolimits_{2})$ have a co-prod\-uct, an
antipode, and a co-unit.

However, there are two ways of choosing the Hopf structure of the
$q$-de\-formed partial derivatives. It is so because the two different
co-pro\-ducts of the $q$-de\-formed partial derivatives are related to the two
versions of Leibniz rules given in Eq.~(\ref{DifKalExtEukQuaDreUnk}) or
Eq.~(\ref{DifKalExtEukQuaDreKon}) of the last subchapter. For a better
insight, we note that you can generalize these Leibniz rules by introducing
so-called L-ma\-tri\-ces $\mathcal{L}_{\partial}$ and $\mathcal{\bar{L}%
}_{\partial}$ ($u\in\mathbb{R}_{q}^{3}$):%
\begin{align}
\partial^{A}\hspace{0.01in}u  &  =(\partial_{(1)}^{A}\triangleright
u)\,\partial_{(2)}^{A}=\partial^{A}\triangleright u\hspace{0.01in}%
+\big ((\mathcal{L}_{\partial}){^{A}}_{B}\triangleright u\big )\partial
^{B},\nonumber\\
\hat{\partial}^{A}\hspace{0.01in}u  &  =(\hat{\partial}_{(\bar{1})}^{A}%
\,\bar{\triangleright}\,u)\,\hat{\partial}_{(\bar{2})}^{A}=\hat{\partial}%
^{A}\,\bar{\triangleright}\,u\hspace{0.01in}+\big ((\mathcal{\bar{L}%
}_{\partial}){^{A}}_{B}\triangleright u\big )\hat{\partial}^{B}.
\label{AllVerRelParAblEle1}%
\end{align}
You can see from the above identities that the two L-ma\-trices determine the
two co-prod\-ucts\footnote{We write the co-product in the so-called Sweedler
notation, i.~e. $\Delta(a)=a_{(1)}\otimes a_{(2)}.$} of the $q$-de\-formed
partial derivatives \cite{ogievetsky1992}:%
\begin{align}
\partial_{(1)}^{A}\otimes\partial_{(2)}^{A}  &  =\partial^{A}\otimes
1+(\mathcal{L}_{\partial}){^{A}}_{B}\otimes\partial^{B},\nonumber\\
\hat{\partial}_{(1)}^{A}\otimes\hat{\partial}_{(2)}^{A}  &  =\hat{\partial
}^{A}\otimes1+(\mathcal{\bar{L}}_{\partial}){^{A}}_{B}\otimes\hat{\partial
}^{B}. \label{KopParAllg}%
\end{align}

The entries of the two L-ma\-tri\-ces consist of generators of the Hopf
algebra $\mathcal{U}_{q}(\operatorname*{su}_{2})$ and powers of a scaling
operator $\Lambda$ [also see Eq.~(\ref{SkaOpeWir})]. For this reason, the
L-ma\-tri\-ces can act on any element of $\mathbb{R}_{q}^{3}$. In this
respect, we say that an element of $\mathbb{R}_{q}^{3}$ has trivial braiding
if the L-ma\-tri\-ces act on it as follows:%
\begin{align}
(\mathcal{L}_{\partial}){^{A}}_{B}\triangleright u  &  =\delta_{B}^{A}%
\hspace{0.01in}u,\nonumber\\
(\mathcal{\bar{L}}_{\partial}){^{A}}_{B}\triangleright u  &  =\delta_{B}%
^{A}\hspace{0.01in}u. \label{trivBrai}%
\end{align}

In Ref.~\cite{Bauer:2003} and Ref.~\cite{Mikulovic:2006}, we have written down
the co-prod\-ucts of the partial derivatives $\partial^{A}$ or $\hat{\partial
}^{A}$, $A\in\{+,3,-\}$, explicitly. By taking into account
Eq.~(\ref{KopParAllg}), you can read off the entries of the L-ma\-tri\-ces
$\mathcal{L}_{\partial}$ and $\mathcal{\bar{L}}_{\partial}$ from these
co-prod\-ucts. In doing so, you find, for example:\footnote{Instead of $T^{3}%
$, one often uses $\tau=1-\lambda T^{3}$.}%
\begin{equation}
(\mathcal{L}_{\partial}){^{-}}_{-}=\Lambda^{1/2}\hspace{0.01in}\tau
^{-1/2}\text{ and }(\mathcal{\bar{L}}_{\partial}){^{+}}_{+}=\Lambda
^{-1/2}\hspace{0.01in}\tau^{-1/2}.
\end{equation}

The scaling operator $\Lambda$ acts on the spatial coordinates or the
corresponding partial derivatives as follows:%
\begin{equation}
\Lambda\triangleright X^{A}=q^{4}X^{A},\qquad\Lambda\triangleright\partial
^{A}=q^{-4}\partial^{A}. \label{SkaOpeWir}%
\end{equation}
These actions imply the commutation relations%
\begin{equation}
\Lambda\hspace{0.01in}X^{A}=q^{4}X^{A}\Lambda,\qquad\Lambda\hspace
{0.01in}\partial^{A}=q^{-4}\partial^{A}\Lambda\label{VerSkaKooQuaEukDrei}%
\end{equation}
if we take into account the Hopf structure of $\Lambda$ \cite{ogievetsky1992}:%
\begin{equation}
\Delta(\Lambda)=\Lambda\otimes\Lambda,\qquad S(\Lambda)=\Lambda^{-1}%
,\qquad\varepsilon(\Lambda)=1. \label{HopStrLam}%
\end{equation}

The Hopf structure of the partial derivatives includes a co-prod\-uct as well
as an antipode and a co-unit. Regarding the co-unit of the partial
derivatives, the following holds \cite{ogievetsky1992}:%
\[
\varepsilon(\partial^{A})=0.
\]
We can obtain the antipodes of the partial derivatives from their
co-prod\-ucts using the following Hopf algebra axioms:%
\begin{equation}
a_{(1)}\cdot S(a_{(2)})=\varepsilon(a)=S(a_{(1)})\cdot a_{(2)}.
\end{equation}
Due to this axiom, we have:%
\begin{equation}
S(\partial^{A})=-\hspace{0.01in}S(\mathcal{L}_{\partial}){^{A}}_{B}%
\hspace{0.01in}\partial^{B},\qquad\bar{S}(\hat{\partial}^{A})=-\hspace
{0.01in}S^{-1}(\mathcal{\bar{L}}_{\partial}){^{A}}_{B}\hspace{0.01in}%
\hat{\partial}^{B}. \label{AlgForInvAntAbl}%
\end{equation}
This way, for example, we get the following expressions for the antipodes of
the partial derivatives $\partial^{-}$ and $\hat{\partial}^{+}$ (also see
Ref.~\cite{Bauer:2003}):%
\begin{equation}
S(\partial^{-})=-\Lambda^{-1/2}\hspace{0.01in}\tau^{1/2}\hspace{0.01in}%
\partial^{-},\qquad\bar{S}(\hat{\partial}^{+})=-\Lambda^{1/2}\hspace
{0.01in}\tau^{1/2}\hspace{0.01in}\hat{\partial}^{+}.
\end{equation}

The antipodes of partial derivatives allow us to write the left actions of
partial derivatives as right actions, and vice versa. Concretely, we have%
\begin{align}
\partial^{A}\triangleright f  &  =f\triangleleft S(\partial^{A}%
)=-\big (f\triangleleft S(\mathcal{L}_{\partial}){^{A}}_{B}\big )\triangleleft
\partial^{B}\nonumber\\
&  =-\big ((\mathcal{L}_{\partial}){^{A}}_{B}\triangleright
f\big )\triangleleft\partial^{B},\nonumber\\[0.06in]
\hat{\partial}^{A}\,\bar{\triangleright}\,f  &  =f\,\bar{\triangleleft}%
\,\bar{S}(\hat{\partial}^{A})=-\big (f\triangleleft S(\mathcal{\bar{L}%
}_{\partial}){^{A}}_{B}\big )\,\bar{\triangleleft}\,\hat{\partial}%
^{B}\nonumber\\
&  =-\big ((\mathcal{\bar{L}}_{\partial}){^{A}}_{B}\triangleright
f\big )\,\bar{\triangleleft}\,\hat{\partial}^{B}, \label{LinkRechtDarN}%
\end{align}
and%
\begin{align}
f\triangleleft\hat{\partial}^{A}  &  =S^{-1}(\hat{\partial}^{A})\triangleright
f=-\hspace{0.01in}\hat{\partial}^{B}\triangleright\big (S^{-1}(\mathcal{L}%
_{\partial}){^{A}}_{B}\triangleright f\big )\nonumber\\
&  =-\hspace{0.01in}\hat{\partial}^{B}\triangleright\big (f\triangleleft
(\mathcal{L}_{\partial}){^{A}}_{B}\big ),\nonumber\\[0.06in]
f\,\bar{\triangleleft}\,\partial^{A}  &  =\bar{S}^{-1}(\partial^{A}%
)\,\bar{\triangleright}\,f=-\hspace{0.01in}\partial^{B}\,\bar{\triangleright
}\,\big (S^{-1}(\mathcal{\bar{L}}_{\partial}){^{A}}_{B}\triangleright
f\big )\nonumber\\
&  =-\hspace{0.01in}\partial^{B}\,\bar{\triangleright}\,\big (f\triangleleft
(\mathcal{\bar{L}}_{\partial}){^{A}}_{B}\big ). \label{RechtsLinksDarN}%
\end{align}

\section{Schr\"{o}dinger equations for a $q$-deformed nonrelativistic
particle\label{KapHerSchGle}}

In Ref.~\cite{Wachter:2020B}, we have chosen the following expression as
Hamilton operator for a free nonrelativistic particle with mass $m$:%
\begin{equation}
H_{0}=-(2\hspace{0.01in}m)^{-1}g_{AB}\hspace{0.01in}\partial^{A}\partial
^{B}=-(2\hspace{0.01in}m)^{-1}\partial^{A}\partial_{A}. \label{Ham2}%
\end{equation}
This choice ensures that $H_{0}$ behaves like a scalar under the action of the
Hopf algebra $\mathcal{U}_{q}(\operatorname*{su}\nolimits_{2})$. Moreover,
$H_{0}$ has trivial braiding if we require:%
\[
\Lambda\hspace{0.01in}m=q^{-8}m\hspace{0.01in}\Lambda.
\]
Due to its definition, the Hamilton operator $H_{0}$ is also a central element
of the algebra of $q$-de\-formed partial derivatives:%
\begin{equation}
\lbrack H_{0},\partial^{A}]=0,\quad A\in\{+,3,-\}. \label{ComHP}%
\end{equation}
The conjugation properties of the partial derivatives imply that $H_{0}$ is
invariant under conjugation [cf. Eq.~(\ref{KonAbl}) of Chap.~\ref{KapParDer}]:%
\begin{equation}
\overline{H_{0}}=H_{0}. \label{RelBedHamFre}%
\end{equation}

We can add a potential $V(\mathbf{x})$ to the free Hamilton operator:%
\begin{equation}
H=H_{0}+V(\mathbf{x}). \label{HamOpePot}%
\end{equation}
Under the action of the Hopf algebra $\mathcal{U}_{q}(\operatorname*{su}%
\nolimits_{2})$, $V(\mathbf{x})$ should behave like a scalar. We also require
$V(\mathbf{x})$ to have trivial braiding and to be central in the algebra of
position space:%
\begin{equation}
V(\mathbf{x})\circledast f(\mathbf{x})=f(\mathbf{x})\circledast V(\mathbf{x}).
\label{CenEle}%
\end{equation}
Moreover, $V(\mathbf{x})$ has to be invariant under conjugation:%
\begin{equation}
\overline{V(\mathbf{x})}=V(\mathbf{x}).
\end{equation}

If we deal with a charged particle moving in the presence of a magnetic field,
we start from the Hamilton operator%
\begin{equation}
H=-(2\hspace{0.01in}m)^{-1}\hspace{0.01in}D^{C}D_{C}, \label{HamVekPot}%
\end{equation}
which depends on the covariant derivatives%
\begin{equation}
D^{C}\hspace{-0.01in}=\partial^{C}\hspace{-0.01in}-\text{i}\hspace
{0.01in}e\hspace{0.01in}A^{C}(\mathbf{x},t). \label{DefKanImp}%
\end{equation}
$A^{C}(\mathbf{x},t)$ denotes the component of the vector potential and $e$ is
the charge of the particle. To ensure that $H$ has well-de\-fined braiding
properties, the L-ma\-tri\-ces have to act on the component $A^{C}%
(\mathbf{x},t)$ in the same way as on the partial derivative $\partial^{C}$.
Since $H$ has to be invariant under quantum space conjugation, we also require%
\begin{equation}
\overline{A^{D}(\mathbf{x},t)}=A_{D}(\mathbf{x},t).
\end{equation}

We can use the Hamilton operators given in Eq.~(\ref{HamOpePot}) or
Eq.~(\ref{HamVekPot}) to write down Schr\"{o}\-dinger equations
\cite{Wachter:2020A}:%
\begin{align}
\text{i}\partial_{t}\triangleright\psi_{R}(\mathbf{x},t)  &  =H\triangleright
\psi_{R}(\mathbf{x},t), & \psi_{L}(\mathbf{x},t)\,\bar{\triangleleft
}\,\partial_{t}\text{i}  &  =\psi_{L}(\mathbf{x},t)\,\bar{\triangleleft
}\,H,\nonumber\\
\text{i}\partial_{t}\,\bar{\triangleright}\,\psi_{R}^{\ast}(\mathbf{x},t)  &
=H\,\bar{\triangleright}\,\psi_{R}^{\ast}(\mathbf{x},t), & \psi_{L}^{\ast
}(\mathbf{x},t)\triangleleft\partial_{t}\text{i}  &  =\psi_{L}^{\ast
}(\mathbf{x},t)\triangleleft H. \label{SchGleQDef1N}%
\end{align}
As was shown in Ref.~\cite{Wachter:2020A}, the time evolution operator for the
quantum space $\mathbb{R}_{q}^{3}$ is of the same form as in the undeformed
case. For this reason, we get solutions to our $q$-de\-formed
Schr\"{o}\-dinger equations by applying the operators $\exp(-$i$tH)$ and
$\exp($i$tH)$ to wave functions at time $t=0$:%
\begin{align}
\psi_{R}(\mathbf{x},t)  &  =\exp(-\text{i}tH)\triangleright\psi_{R}%
(\mathbf{x},0), & \psi_{L}(\mathbf{x},t)  &  =\psi_{L}(\mathbf{x}%
,0)\,\bar{\triangleleft}\,\exp(\text{i}Ht),\nonumber\\
\psi_{R}^{\ast}(\mathbf{x},t)  &  =\exp(-\text{i}tH)\,\bar{\triangleright
}\,\psi_{R}^{\ast}(\mathbf{x},0), & \psi_{L}^{\ast}(\mathbf{x},t)  &
=\psi_{L}^{\ast}(\mathbf{x},0)\triangleleft\exp(\text{i}Ht).
\label{SchGleHamOpe1}%
\end{align}

Last but not least, we require that the solutions to the Schr\"{o}\-dinger
equations in Eq.~(\ref{SchGleQDef1N}) behave as follows under quantum space
conjugation:%
\begin{equation}
\overline{\psi_{L}(\mathbf{x},t)}=\psi_{R}(\mathbf{x},t),\text{\qquad
}\overline{\psi_{L}^{\ast}(\mathbf{x},t)}=\psi_{R}^{\ast}(\mathbf{x},t).
\label{KonEigSchWelWdh}%
\end{equation}
This condition ensures that quantum space conjugation transforms the
Schr\"{o}\-dinger equations on the left side of Eq.~(\ref{SchGleQDef1N}) into
the Schr\"{o}\-dinger equations on the right side of Eq.~(\ref{SchGleQDef1N})
and vice versa.

\section{$q$-Versions of Green's theorem\label{KapGre}}

We need $q$-versions of Green's theorem to derive
conservation laws. We will show in this chapter how to get these
$q$-de\-formed analogs from the Leibniz rules for $q$-de\-formed partial
derivatives and the properties of L-ma\-tri\-ces.\footnote{A one-dimensional $q$-version of Green's theorem was
given in Ref.~\cite{Cerchiai:1999}.}

The Leibniz rules in Eq.~(\ref{AllVerRelParAblEle1}) of Chap.~\ref{KapParDer}
imply the following rule for left actions of $q$-de\-formed partial
derivatives:%
\begin{align}
\psi\circledast\partial^{A}\triangleright\phi &  =\partial_{(2)}%
^{A}\triangleright\big [\psi\triangleleft\partial_{(1)}^{A}\circledast
\phi\big ]\nonumber\\
&  =\partial^{B}\triangleright\left[  \psi\triangleleft(\mathcal{L}_{\partial
}){^{A}}_{\hspace{-0.01in}B}\circledast\phi\right]  +\psi\triangleleft
\partial^{A}\circledast\phi.\label{UmgLeiReg1}%
\end{align}
We can prove the identity above in the following way:%
\begin{align}
&  \partial^{B}\triangleright\left[  \psi\triangleleft(\mathcal{L}_{\partial
}){^{A}}_{\hspace{-0.01in}B}\circledast\phi\right]  =\partial^{B}%
\triangleright\left[  S^{-1}(\mathcal{L}_{\partial}){^{A}}_{\hspace{-0.01in}%
B}\triangleright\psi\circledast\phi\right]  \nonumber\\
&  \qquad\qquad=\partial^{B}S^{-1}(\mathcal{L}_{\partial}){^{A}}%
_{\hspace{-0.01in}B}\triangleright\psi\circledast\phi+(\mathcal{L}_{\partial
}){^{B}}_{\hspace{-0.01in}C}\hspace{0.02in}S^{-1}(\mathcal{L}_{\partial}%
){^{A}}_{\hspace{-0.01in}B}\triangleright\psi\circledast\partial
^{C}\triangleright\phi\nonumber\\
&  \qquad\qquad=-\hspace{0.01in}\psi\triangleleft\partial^{A}\circledast
\phi+\delta_{C}^{A}\hspace{0.02in}\psi\circledast\partial^{C}\triangleright
\phi\nonumber\\
&  \qquad\qquad=-\hspace{0.01in}\psi\triangleleft\partial^{A}\circledast
\phi+\psi\circledast\partial^{A}\triangleright\phi.
\end{align}
In the first step, we have used the fact that the entries $(\mathcal{L}%
_{\partial}){^{A}}_{\hspace{-0.01in}B}$ are elements of a Hopf algebra (also
see Chap.~\ref{KapHofStr}). Thus, we can write the right action as a
left-action by using the inverse antipode of this Hopf algebra. The second
step results from Eq.~(\ref{AllVerRelParAblEle1}) of Chap.~\ref{KapParDer}.
The third step and the fourth step are a consequence of the Hopf algebra axiom%
\begin{equation}
a_{(2)}\cdot S^{-1}a_{(1)}=\varepsilon(a)
\end{equation}
and the following properties of L-ma\-tri\-ces:%
\begin{equation}
\Delta(\mathcal{L}_{\partial}){^{A}}_{\hspace{-0.01in}B}=(\mathcal{L}%
_{\partial}){^{A}}_{\hspace{-0.01in}C}\otimes(\mathcal{L}_{\partial}){^{C}%
}_{\hspace{-0.01in}B},\qquad\varepsilon(\mathcal{L}_{\partial}){^{A}}%
_{\hspace{-0.01in}B}=\delta_{B}^{A}.\label{HopLMat}%
\end{equation}
Similar arguments lead to%
\begin{align}
\psi\triangleleft\partial^{A}\circledast\phi &  =[\psi\circledast
\partial_{(1)}^{A}\triangleright\phi]\triangleleft\partial_{(2)}%
^{A}\nonumber\\
&  =[\psi\circledast(\mathcal{L}_{\partial}){^{A}}_{\hspace{-0.01in}%
B}\triangleright\phi]\triangleleft\partial^{B}+\psi\circledast\partial
^{A}\triangleright\phi.\label{UmgLeiReg2}%
\end{align}

We can now derive a $q$-version of Green's theorem. Due to
Eq.~(\ref{UmgLeiReg1}), we have%
\begin{equation}
\psi\triangleleft\partial^{A}\partial_{A}\circledast\phi=\psi\triangleleft
\partial^{A}\circledast\partial_{A}\triangleright\phi-\partial^{C}%
\triangleright\lbrack\psi\triangleleft\partial^{A}(\mathcal{L}_{\partial
}){^{B}}_{\hspace{-0.01in}C}\circledast\phi]\hspace{0.02in}g_{AB}
\label{ZwiRecWskErhDreQua1}%
\end{equation}
and%
\begin{equation}
\psi\circledast\partial^{A}\partial_{A}\triangleright\phi=\psi\triangleleft
\partial^{A}\circledast\partial_{A}\triangleright\phi+\partial^{C}%
\triangleright\lbrack\psi\triangleleft(\mathcal{L}_{\partial}){^{A}}%
_{\hspace{-0.01in}C}\circledast\partial^{\hspace{0.01in}B}\triangleright
\phi]\hspace{0.02in}g_{AB}. \label{ZwiRecWskErhDreQua2}%
\end{equation}
Combining Eq.~(\ref{ZwiRecWskErhDreQua1}) and Eq.~(\ref{ZwiRecWskErhDreQua2}),
we obtain:%
\begin{align}
&  \psi\triangleleft\partial^{A}\partial_{A}\circledast\phi-\psi
\circledast\partial^{A}\partial_{A}\triangleright\phi=\nonumber\\
&  \qquad=-\hspace{0.01in}\partial^{C}\triangleright\left[  \psi
\triangleleft\partial^{A}(\mathcal{L}_{\partial}){^{B}}_{\hspace{-0.01in}%
C}\circledast\phi+\psi\triangleleft(\mathcal{L}_{\partial}){^{A}}%
_{\hspace{-0.01in}C}\circledast\partial^{B}\triangleright\phi\right]
\hspace{0.02in}g_{AB}\nonumber\\
&  \qquad=-\hspace{0.01in}\partial^{C}\triangleright\left[  q^{-2}%
\hspace{0.01in}\psi\triangleleft(\mathcal{L}_{\partial}){^{A}}_{\hspace
{-0.01in}C}\hspace{0.02in}\partial^{B}\hspace{-0.01in}\circledast\phi
+\psi\triangleleft(\mathcal{L}_{\partial}){^{A}}_{\hspace{-0.01in}%
C}\circledast\partial^{B}\triangleright\phi\right]  \hspace{0.02in}g_{AB}.
\label{QGreIde0}%
\end{align}
The last step of the above calculation results from the following identities:%
\begin{align}
g_{AB}\hspace{0.01in}\partial^{A}(\mathcal{L}_{\partial}){^{B}}_{\hspace
{-0.01in}C}  &  =g_{AB}\hspace{0.01in}(\mathcal{L}_{\partial}){^{D}}%
_{\hspace{-0.01in}C}\hspace{0.01in}[\hspace{0.01in}\partial^{A}\triangleleft
(\mathcal{L}_{\partial}){^{B}}_{\hspace{-0.01in}D}]\nonumber\\
&  =g_{AB}\hspace{0.01in}(\mathcal{L}_{\partial}){^{D}}_{\hspace{-0.01in}%
C}\hspace{0.01in}q^{4}\hat{R}{}{^{\hspace{0.01in}AB}}_{\hspace{-0.01in}%
DE}\,\partial^{E}=q^{-6}q^{4}g_{DE}\hspace{0.01in}(\mathcal{L}_{\partial
}){^{D}}_{\hspace{-0.01in}C}\hspace{0.01in}\partial^{E}\nonumber\\
&  =q^{-2}(\mathcal{L}_{\partial}){^{A}}_{\hspace{-0.01in}C}\hspace
{0.02in}\partial^{B}g_{AB}\hspace{0.01in}. \label{VerAplLMat}%
\end{align}
The first identity of Eq.~(\ref{VerAplLMat}) follows from the co-prod\-uct of
the L-ma\-tri\-ces [cf. Eq.~(\ref{HopLMat})] and the commutation relations of
the Hopf algebra $\mathbb{R}_{q}^{3}\rtimes\hspace{0.01in}\mathcal{U}%
_{q}(\operatorname*{su}\nolimits_{2})$:%
\begin{equation}
\partial\cdot u=u_{(2)}\cdot(\partial\triangleleft u_{(1)}).
\end{equation}
The second identity of Eq.~(\ref{VerAplLMat}) holds because the vector
representation of the R-ma\-trix determines the action of $(\mathcal{L}%
_{\partial}){^{A}}_{\hspace{-0.01in}B}$ on a vector. The penultimate identity
of Eq.~(\ref{VerAplLMat}) follows from the projector decomposition of the
R-ma\-trix of the $q$-de\-formed Euclidean space \cite{Lorek:1997eh},
i.~e.\footnote{The projector $P_{A}$ is a $q$-analog of an antisymmetrizer,
the projector $P_{S}$ is the $q$-deformed trace-free symmetrizer, and $P_{T}$
is the $q$-de\-formed trace-projector.}%
\begin{equation}
\hat{R}=P_{S}+q^{-6}P_{T}-q^{-4}P_{A},
\end{equation}
if we take into account:%
\begin{equation}
g_{AB}\hspace{0.01in}(P_{S}){^{AB}}_{CD}=g_{AB}\hspace{0.01in}(P_{A}){^{AB}%
}_{CD}=0,\qquad g_{AB}\hspace{0.01in}(P_{T}){^{AB}}_{CD}=g_{CD}.
\end{equation}
By similar reasonings using Eq.~(\ref{UmgLeiReg2}), we can also derive the
following identity:%
\begin{align}
&  \psi\triangleleft\partial^{A}\partial_{A}\circledast\phi-\psi
\circledast\partial^{A}\partial_{A}\triangleright\phi=\nonumber\\
&  \qquad=g_{AB}\hspace{-0.01in}\left[  \psi\triangleleft\partial
^{A}\circledast(\mathcal{L}_{\partial}){^{B}}_{\hspace{-0.01in}C}%
\triangleright\phi+q^{2}\hspace{0.01in}\psi\circledast\partial^{A}%
(\mathcal{L}_{\partial}){^{B}}_{\hspace{-0.01in}C}\triangleright\phi\right]
\triangleleft\partial^{C}. \label{QGreIde1}%
\end{align}

If there is a vector potential, we have to substitute the partial derivatives
$\partial^{C}$ by the operators $D^{C}$ [cf. Eq.~(\ref{DefKanImp}) of
Chap.~\ref{KapHerSchGle}]. Instead of Eq.~(\ref{QGreIde0}) and
Eq.~(\ref{QGreIde1}), we have the identities
\begin{align}
&  \psi\triangleleft D^{C}D_{C}\circledast\phi-\psi\circledast D^{C}%
D_{C}\triangleright\phi=\nonumber\\
&  \qquad=-\hspace{0.02in}\partial^{F}\triangleright\left[  q^{-2}%
\hspace{0.01in}\psi\triangleleft(\mathcal{L}_{\partial}){^{B}}_{\hspace
{-0.01in}F}\hspace{0.01in}D^{C}\circledast\phi+\psi\triangleleft
(\mathcal{L}_{\partial}){^{B}}_{\hspace{-0.01in}F}\circledast D^{C}%
\triangleright\phi\right]  \hspace{0.02in}g_{BC}\nonumber\\
&  \qquad=g_{BC}\hspace{-0.01in}\left[  \psi\triangleleft D^{B}\circledast
(\mathcal{L}_{\partial}){^{C}}_{\hspace{-0.01in}F}\triangleright\phi
+q^{2}\hspace{0.01in}\psi\circledast D^{B}(\mathcal{L}_{\partial}){^{C}%
}_{\hspace{-0.01in}F}\triangleright\phi\right]  \triangleleft\partial^{F}
\label{GreIdeKinImp}%
\end{align}
with%
\begin{align}
D^{C}\triangleright\phi &  =\partial^{C}\triangleright\phi-\text{i}%
\hspace{0.01in}eA^{C}\circledast\phi,\nonumber\\
\psi\triangleleft D^{C}  &  =\psi\triangleleft\partial^{C}-\text{i}%
\hspace{0.01in}\psi\circledast A^{C}e.
\end{align}
We can derive these identities in the same way as the $q$-versions of Green's
theorem if we use the formula%
\begin{align}
\psi\circledast D^{B}\triangleright\phi-\psi\triangleleft D^{B}\circledast\phi
&  =\partial^{C}\triangleright\left[  \psi\triangleleft(\mathcal{L}_{\partial
}){^{B}}_{\hspace{-0.01in}C}\circledast\phi\right] \nonumber\\
&  =-\left[  \psi\circledast(\mathcal{L}_{\partial}){^{B}}_{\hspace{-0.01in}%
C}\triangleright\phi\right]  \triangleleft\partial^{C} \label{UmgLeiRegKinImp}%
\end{align}
and recall that the rules in\ Eq.~(\ref{AllVerRelParAblEle1}) of
Chap.~\ref{KapHofStr}\ also hold for the operators $D^{C}$.

\section{Conservation of probability\label{KonGleWskErhKap}}

We require that the solutions to the free $q$-de\-formed Schr\"{o}\-dinger
equations given in Eq.~(\ref{SchGleQDef1N}) of Chap.~\ref{KapHerSchGle} are
subject to the following normalization condition:%
\begin{equation}
1=\frac{1}{2}\int\text{d}_{q}^{3}\hspace{0.01in}x\left(  \psi_{L}^{\ast
}(\mathbf{x},t)\circledast\psi_{R}(\mathbf{x},t)+\psi_{L}(\mathbf{x}%
,t)\circledast\psi_{R}^{\ast}(\mathbf{x},t)\right)  . \label{NorBed}%
\end{equation}
This condition confirms that the probability density for a nonrelativistic
particle depends on the following expressions:%
\begin{align}
\rho(\mathbf{x},t)  &  =\psi_{L}^{\ast}(\mathbf{x},t)\circledast\psi
_{R}(\mathbf{x},t),\nonumber\\
\rho^{\ast}(\mathbf{x},t)  &  =\psi_{L}(\mathbf{x},t)\circledast\psi_{R}%
^{\ast}(\mathbf{x},t). \label{DefProp4Hab}%
\end{align}
Due to Eq.~(\ref{KonEigSchWelWdh}) of Chap.~\ref{KapHerSchGle}, the
expressions above transform into each other by conjugation:%
\begin{equation}
\overline{\rho(\mathbf{x},t)}=\rho^{\ast}(\mathbf{x},t). \label{ConProbDen}%
\end{equation}

Since the wave functions stay normalized as they evolve in time,
$\rho(\mathbf{x},t)$ and $\rho^{\ast}(\mathbf{x},t)$ must satisfy continuity
equations. In the following, we will derive these continuity
equations.\footnote{The considerations of the present chapter and the
following ones are similar to those in Ref.~\cite{Pauli:1980}.} We first
consider a nonrelativistic particle in an external force field with a scalar
potential $V(\mathbf{x})$. In this case, we use the Hamilton operator $H$
given by Eq.~(\ref{HamOpePot}) of Chap.~\ref{KapHerSchGle}. We calculate the
time derivative of $\rho(\mathbf{x},t)$ by taking into account the
Schr\"{o}\-dinger equations in Eq.~(\ref{SchGleQDef1N}) of Chap.
\ref{KapHerSchGle}:%
\begin{align}
\partial_{t}\hspace{0.01in}\triangleright\rho(\mathbf{x},t)  &  =\partial
_{t}\triangleright\psi_{L}^{\ast}(\mathbf{x},t)\circledast\psi_{R}%
(\mathbf{x},t)+\psi_{L}^{\ast}(\mathbf{x},t)\circledast\partial_{t}%
\triangleright\psi_{R}(\mathbf{x},t)\nonumber\\
&  =\text{i}\big (\psi_{L}^{\ast}\triangleleft H\circledast\psi_{R}-\psi
_{L}^{\ast}\circledast H\triangleright\psi_{R}\big )\nonumber\\
&  =\text{i}\big (\psi_{L}^{\ast}\triangleleft H_{0}\circledast\psi_{R}%
-\psi_{L}^{\ast}\circledast H_{0}\triangleright\psi_{R}\big )\nonumber\\
&  =-\hspace{0.01in}\text{i}\big (\psi_{L}^{\ast}\triangleleft\partial
^{A}\partial_{A}(2\hspace{0.01in}m)^{-1}\hspace{-0.01in}\circledast\psi
_{R}-\psi_{L}^{\ast}\circledast(2\hspace{0.01in}m)^{-1}\partial^{A}%
\partial_{A}\triangleright\psi_{R}\big ). \label{UmfKonGleEuk1}%
\end{align}
Note that the contributions due to $V(\mathbf{x})$ cancel each other out.
Applying one of the $q$-versions of Green's theorem [cf.
Eq.~(\ref{QGreIde0}) of Chap.~\ref{KapGre}] to the last expression in
Eq.~(\ref{UmfKonGleEuk1}), we obtain a \textit{continuity equation for the
probability density} $\rho(\mathbf{x},t)$, i.~e.%
\begin{equation}
\partial_{t}\hspace{0.01in}\triangleright\rho(\mathbf{x},t)+\partial
^{A}\triangleright j_{A}(\mathbf{x},t)=0 \label{KonGleDreDim1}%
\end{equation}
with the \textit{probability current}%
\begin{align}
j_{A}(\mathbf{x},t)=  &  -\frac{\text{i}}{2\hspace{0.01in}m}\hspace
{0.01in}q^{-2}\hspace{0.01in}\psi_{L}^{\ast}(\mathbf{x},t)\triangleleft
(\mathcal{L}_{\partial}){^{B}}_{\hspace{-0.01in}A}\,\partial_{B}%
\hspace{0.01in}\circledast\psi_{R}(\mathbf{x},t)\nonumber\\
&  -\frac{\text{i}}{2\hspace{0.01in}m}\hspace{0.01in}\psi_{L}^{\ast
}(\mathbf{x},t)\triangleleft(\mathcal{L}_{\partial}){^{B}}_{\hspace{-0.01in}%
A}\circledast\partial_{B}\triangleright\psi_{R}(\mathbf{x},t).
\label{StrDicVekKov}%
\end{align}

The continuity equation for the probability density $\rho^{\ast}%
(\mathbf{x},t)$ can be derived by similar considerations as above or by
conjugating Eq.~(\ref{KonGleDreDim1}). This way, we get%
\begin{equation}
\rho^{\ast}(\mathbf{x},t)\,\bar{\triangleleft}\,\partial_{t}+(j^{\ast}%
)^{A}(\mathbf{x},t)\,\bar{\triangleleft}\,\partial_{A}=0
\label{KonGleDreDim1Ste}%
\end{equation}
with the new probability current%
\begin{align}
(j^{\ast})^{A}(\mathbf{x},t)=  &  -\frac{\text{i}}{2\hspace{0.01in}m}%
\hspace{0.01in}q^{-2}g_{BC}\hspace{0.01in}\psi_{L}(\mathbf{x},t)\circledast
\partial^{B}\hspace{-0.01in}(\mathcal{\bar{L}}_{\partial}){^{C}}%
_{\hspace{-0.01in}F}\,\bar{\triangleright}\,\psi_{R}^{\ast}(\mathbf{x}%
,t)\hspace{0.01in}g^{FA}\nonumber\\
&  -\frac{\text{i}}{2\hspace{0.01in}m}\hspace{0.01in}g_{BC}\hspace{0.01in}%
\psi_{L}(\mathbf{x},t)\,\bar{\triangleleft}\,\partial^{B}\circledast
(\mathcal{\bar{L}}_{\partial}){^{C}}_{\hspace{-0.01in}F}\,\bar{\triangleright
}\,\psi_{R}^{\ast}(\mathbf{x},t)\hspace{0.01in}g^{FA}. \label{StrDicVekKovSte}%
\end{align}
In analogy to Eq.~(\ref{ConProbDen}), it holds:%
\begin{equation}
\overline{j_{A}(\mathbf{x},t)}=(j^{\ast})^{A}(\mathbf{x},t).
\label{ConProbCur}%
\end{equation}

Next, we derive the continuity equation for the probability density of a
charged particle in a magnetic field. We calculate the time derivative of the
probability density by taking into account the Hamilton operator in
Eq.~(\ref{HamVekPot}) of Chap.~\ref{KapHerSchGle} and the identity in
Eq.~(\ref{GreIdeKinImp}) of Chap.~\ref{KapGre}:%
\begin{align}
\partial_{t}\hspace{0.01in}\triangleright\rho(\mathbf{x},t)=  &
\,-\text{i}\hspace{0.01in}\psi_{L}^{\ast}\triangleleft D^{C}D_{C}%
\hspace{-0.01in}\circledast\psi_{R}+\text{i}\hspace{0.01in}\psi_{L}^{\ast
}\circledast(2\hspace{0.01in}m)^{-1}D^{C}D_{C}\triangleright\psi
_{R}\nonumber\\
=  &  -\hspace{0.02in}\text{i}^{-1}(2\hspace{0.01in}m)^{-1}\partial
^{F}\triangleright\left[  q^{-2}\hspace{0.01in}\psi_{L}^{\ast}\triangleleft
(\mathcal{L}_{\partial}){^{B}}_{\hspace{-0.01in}F}\hspace{0.01in}%
D^{C}\circledast\psi_{R}\right]  \hspace{0.02in}g_{BC}\nonumber\\
&  -\hspace{0.02in}\text{i}^{-1}(2\hspace{0.01in}m)^{-1}\partial
^{F}\triangleright\left[  \psi_{L}^{\ast}\triangleleft(\mathcal{L}_{\partial
}){^{B}}_{\hspace{-0.01in}F}\circledast D^{C}\triangleright\psi_{R}\right]
\hspace{0.02in}g_{BC}.
\end{align}
From the last expression, we can read off the probability current for a
charged particle in a magnetic field:%
\begin{align}
j_{A}(\mathbf{x},t)=  &  -\frac{\text{i}}{2\hspace{0.01in}m}\hspace
{0.01in}q^{-2}\hspace{0.01in}\psi_{L}^{\ast}(\mathbf{x},t)\triangleleft
(\mathcal{L}_{\partial}){^{B}}_{\hspace{-0.01in}A}\,D_{B}\hspace
{0.01in}\circledast\psi_{R}(\mathbf{x},t)\nonumber\\
&  -\frac{\text{i}}{2\hspace{0.01in}m}\hspace{0.01in}\psi_{L}^{\ast
}(\mathbf{x},t)\triangleleft(\mathcal{L}_{\partial}){^{B}}_{\hspace{-0.01in}%
A}\circledast D_{B}\triangleright\psi_{R}(\mathbf{x},t).
\end{align}
We can also obtain the expression above from that given in
Eq.~(\ref{StrDicVekKov}) if we replace $\partial^{C}$ by $D^{C}$. Using the
explicit form of $D^{C}$ [cf. Eq.~(\ref{DefKanImp}) of
Chap.~\ref{KapHerSchGle}], the probability current of a charged particle in a
magnetic field can be witten as follows:%
\begin{align}
j_{A}=  &  -\frac{\text{i}q^{-2}}{2\hspace{0.01in}m}\,q^{-2}\hspace
{0.01in}\psi_{L}^{\ast}\triangleleft(\mathcal{L}_{\partial}){^{B}}%
_{\hspace{-0.01in}A}\hspace{0.01in}\partial_{B}\circledast\psi_{R}\nonumber\\
&  -\frac{\text{i}}{2\hspace{0.01in}m}\,\psi_{L}^{\ast}\triangleleft
(\mathcal{L}_{\partial}){^{B}}_{\hspace{-0.01in}A}\circledast\partial
_{B}\triangleright\psi_{R}\nonumber\\
&  +\frac{e}{2\hspace{0.01in}m}(q^{-2}+1)\,\psi_{L}^{\ast}\triangleleft
(\mathcal{L}_{\partial}){^{B}}_{\hspace{-0.01in}A}\circledast A_{B}%
\circledast\psi_{R}. \label{StrDicVekKovVekPot}%
\end{align}
There is another expression for the probability current of a particle in a
magnetic field. In analogy to Eq.~(\ref{ConProbCur}), we can get this
expression by conjugation:%
\begin{align}
(j^{\ast})^{A}(\mathbf{x},t)=  &  -\frac{\text{i}}{2\hspace{0.01in}m}%
\hspace{0.01in}q^{-2}g_{BC}\hspace{0.01in}\psi_{L}(\mathbf{x},t)\circledast
\partial^{B}\hspace{-0.01in}(\mathcal{\bar{L}}_{\partial}){^{C}}%
_{\hspace{-0.01in}F}\,\bar{\triangleright}\,\psi_{R}^{\ast}(\mathbf{x}%
,t)\hspace{0.01in}g^{FA}\nonumber\\
&  -\frac{\text{i}}{2\hspace{0.01in}m}\hspace{0.01in}g_{BC}\hspace{0.01in}%
\psi_{L}(\mathbf{x},t)\,\bar{\triangleleft}\,\partial^{B}\circledast
(\mathcal{\bar{L}}_{\partial}){^{C}}_{\hspace{-0.01in}F}\,\bar{\triangleright
}\,\psi_{R}^{\ast}(\mathbf{x},t)\hspace{0.01in}g^{FA}\nonumber\\
&  +\frac{e}{2\hspace{0.01in}m}(q^{-2}+1)\hspace{0.01in}g_{BC}\hspace
{0.01in}\psi_{L}\circledast A^{B}\circledast(\mathcal{\bar{L}}_{\partial
}){^{C}}_{\hspace{-0.01in}F}\,\bar{\triangleright}\,\psi_{R}^{\ast}%
\hspace{0.02in}g^{FA}.
\end{align}

By integrating the above continuity equations, we can show that the wave
functions stay normalized as they evolve in time. The following calculation
using Eq.~(\ref{KonGleDreDim1}) shall serve as an example:%
\begin{align}
\partial_{t}\hspace{0.01in}\triangleright\hspace{-0.01in}\int\text{d}_{q}%
^{3}\hspace{0.01in}x\,\psi_{L}^{\ast}(\mathbf{x},t)\circledast\psi
_{R}(\mathbf{x},t)  &  =\partial_{t}\hspace{0.01in}\triangleright
\hspace{-0.01in}\int\text{d}_{q}^{3}\hspace{0.01in}x\,\rho(\mathbf{x}%
,t)\nonumber\\
&  =-\hspace{-0.01in}\int\text{d}_{q}^{3}\hspace{0.01in}x\,\partial
^{A}\triangleright j_{A}(\mathbf{x},t)=0. \label{ErhWskSchTei}%
\end{align}
The last step of the above calculation follows from the fact that the action
of the derivative $\partial^{A}$ leads to surface terms vanishing at infinity
[cf. Eq.~(\ref{StoThe}) of Chap.~\ref{KapParDer}].

\section{Conservation of momentum\label{KapEneImpErh}}

We introduce the following expressions for \textit{momentum density}:%
\begin{align}
i^{A}  &  =\frac{1}{2\text{i}}\left(  \psi_{L}^{\ast}\circledast\partial
^{A}\triangleright\psi_{R}+\psi_{L}^{\ast}\triangleleft\partial^{A}%
\circledast\psi_{R}\right)  ,\nonumber\\
i_{A}^{\ast}  &  =\frac{1}{2\text{i}}\left(  \psi_{L}\circledast\partial
_{A}\,\bar{\triangleright}\,\psi_{R}^{\ast}+\psi_{L}\,\bar{\triangleleft
}\,\partial_{A}\circledast\psi_{R}^{\ast}\right)  . \label{FreImpDicDefN}%
\end{align}
If we integrate the above expressions over all space and take into account
Eq.~(\ref{PatIntUneRaumInt}) of\ Chap.~\ref{KapParDer}, we get expectation
values for the momentum components of a particle:%
\begin{align}
\langle\hspace{0.01in}p^{A}\rangle &  =\int\text{d}_{q}^{3}\hspace
{0.01in}x\,i^{A}=\int\text{d}_{q}^{3}\hspace{0.01in}x\,\psi_{L}^{\ast
}\circledast\text{i}^{-1}\partial^{A}\triangleright\psi_{R},\nonumber\\
\langle\hspace{0.01in}p_{A}^{\ast}\rangle &  =\int\text{d}_{q}^{3}%
\hspace{0.01in}x\,i_{A}^{\ast}=\int\text{d}_{q}^{3}\hspace{0.01in}x\,\psi
_{L}\circledast\text{i}^{-1}\partial_{A}\,\bar{\triangleright}\,\psi_{R}%
^{\ast}. \label{ExpImp}%
\end{align}
Using Eq.~(\ref{KonEigSchWelWdh}) in Chap.~\ref{KapHerSchGle} together with
Eqs.~(\ref{RegConAbl}), (\ref{PatIntUneRaumInt}), and (\ref{KonEigVolInt}) in
Chap.~\ref{KapParDer}, we can show that quantum space conjugation transforms
the first expression in Eq.~(\ref{FreImpDicDefN}) or Eq.~(\ref{ExpImp}) into
the second one and vice versa:%
\begin{equation}
\overline{i^{A}}=i_{A}^{\ast},\qquad\overline{\langle\hspace{0.01in}%
p^{A}\rangle}=\langle\hspace{0.01in}p_{A}^{\ast}\rangle.
\end{equation}

Next, we calculate the time derivative of the momentum density $i^{A}$. To
this end, we need the Schr\"{o}\-dinger equations given in
Eq.~(\ref{SchGleQDef1N}) of Chap.~\ref{KapHerSchGle}:%
\begin{align}
\partial_{t}\triangleright i^{A}=  &  \,\frac{1}{2}\hspace{0.01in}%
\big (\psi_{L}^{\ast}\triangleleft H\circledast\partial^{A}\triangleright
\psi_{R}-\psi_{L}^{\ast}\circledast\partial^{A}\triangleright(H\triangleright
\psi_{R})\big )\nonumber\\
&  +\frac{1}{2}\hspace{0.01in}\big ((\psi_{L}^{\ast}\triangleleft
H)\triangleleft\partial^{A}\circledast\psi_{R}-\psi_{L}^{\ast}\triangleleft
\partial^{A}\circledast H\triangleright\psi_{R}\big ).
\end{align}
With the explicit form of the Hamilton operator $H$ [cf. Eq.~(\ref{HamOpePot}%
)\ of Chap.~\ref{KapHerSchGle}], we get the following expression for time
derivative $i^{A}$:%
\begin{align}
\partial_{t}\triangleright i^{A}=  &  -\frac{1}{4\hspace{0.01in}m}%
\hspace{0.01in}\big (\psi_{L}^{\ast}\triangleleft\partial^{\hspace{0.01in}%
B}\partial_{B}\circledast\partial^{A}\triangleright\psi_{R}-\psi_{L}^{\ast
}\circledast\partial^{\hspace{0.01in}B}\partial_{B}\hspace{0.01in}\partial
^{A}\triangleright\psi_{R}\big )\nonumber\\
&  -\frac{1}{4\hspace{0.01in}m}\hspace{0.01in}\big (\psi_{L}^{\ast
}\triangleleft\partial^{A}\partial^{\hspace{0.01in}B}\partial_{B}%
\circledast\psi_{R}-\psi_{L}^{\ast}\triangleleft\partial^{A}\circledast
\partial^{\hspace{0.01in}B}\partial_{B}\triangleright\psi_{R}\big )\nonumber\\
&  +\frac{1}{2}\hspace{0.01in}\big (\psi_{L}^{\ast}\circledast V\circledast
\partial^{A}\triangleright\psi_{R}-\psi_{L}^{\ast}\circledast\partial
^{A}\triangleright(V\circledast\psi_{R})\big )\nonumber\\
&  +\frac{1}{2}\hspace{0.01in}\big ((\psi_{L}^{\ast}\circledast
V)\triangleleft\partial^{A}\circledast\psi_{R}-\psi_{L}^{\ast}\triangleleft
\partial^{A}\circledast V\circledast\psi_{R}\big ). \label{ZeiAblImpDic3Dim}%
\end{align}
Since $V$ is invariant under the action of $\mathcal{U}_{q}(\operatorname*{su}%
\nolimits_{2})$ as well as that of $\Lambda$, we have:%
\begin{equation}
(\mathcal{L}_{\partial}){^{A}}_{\hspace{-0.01in}B}\triangleright
V=V\triangleleft(\mathcal{L}_{\partial}){^{A}}_{\hspace{-0.01in}B}%
=\varepsilon((\mathcal{L}_{\partial}){^{A}}_{\hspace{-0.01in}B})\hspace
{0.02in}V=\delta_{B}^{A}\hspace{0.02in}V. \label{TrivBraiV}%
\end{equation}
For this reason, the Leibniz rules of $q$-de\-formed partial derivatives imply
[see Eq.~(\ref{AllVerRelParAblEle1}) of Chap.~\ref{KapParDer}]:%
\begin{align}
\psi_{L}^{\ast}\circledast\partial^{A}\triangleright(V\circledast\psi_{R})  &
=\psi_{L}^{\ast}\circledast(\partial^{A}\triangleright V)\circledast\psi
_{R}+\psi_{L}^{\ast}\circledast V\circledast\partial^{A}\triangleright\psi
_{R},\nonumber\\
(\psi_{L}^{\ast}\circledast V)\triangleleft\partial^{A}\circledast\psi_{R}  &
=\psi_{L}^{\ast}\circledast(V\triangleleft\partial^{A})\circledast\psi
_{R}+\psi_{L}^{\ast}\triangleleft\partial^{A}\circledast V\circledast\psi_{R}.
\end{align}
With these identities, we can combine the last two expressions on the
left-hand side of Eq.~(\ref{ZeiAblImpDic3Dim}) to a \textit{force density}:%
\begin{equation}
f^{A}=-\hspace{0.01in}\frac{1}{2}\big (\psi_{L}^{\ast}\circledast(\partial
^{A}\triangleright V)\circledast\psi_{R}-\psi_{L}^{\ast}\circledast
(V\triangleleft\partial^{A})\circledast\psi_{R}\big ). \label{KraDicDreDim}%
\end{equation}
Moreover, we can write the first two expressions on the left-hand side of
Eq.~(\ref{ZeiAblImpDic3Dim}) as divergence. To achieve this, we apply
Eq.~(\ref{QGreIde0}) of Chap.~\ref{KapGre}\ with%
\begin{equation}
\psi=\psi_{L}^{\ast},\text{\qquad}\phi=\partial^{A}\triangleright\psi_{R}%
\end{equation}
or with%
\begin{equation}
\psi=\psi_{L}^{\ast}\triangleleft\partial^{A},\text{\qquad}\phi=\psi_{R}.
\end{equation}
This way, we get the continuity equation%
\begin{equation}
\partial_{t}\triangleright i_{A}=-\hspace{0.01in}\partial^{\hspace{0.01in}%
B}\triangleright T_{BA}+f_{A} \label{ImpDichAblSpaTen}%
\end{equation}
with the following \textit{stress tensor}:%
\begin{align}
T_{BA}=  &  -\frac{1}{4\hspace{0.01in}m}\hspace{0.01in}q^{-2}\psi_{L}^{\ast
}\triangleleft(\mathcal{L}_{\partial}){^{C}}_{\hspace{-0.01in}B}%
\hspace{0.02in}\partial_{C}\circledast\partial_{A}\triangleright\psi
_{R}\nonumber\\
&  -\frac{1}{4\hspace{0.01in}m}\hspace{0.01in}\psi_{L}^{\ast}\triangleleft
(\mathcal{L}_{\partial}){^{C}}_{\hspace{-0.01in}B}\circledast\partial
_{C}\hspace{0.01in}\partial_{A}\triangleright\psi_{R}\nonumber\\
&  -\frac{1}{4\hspace{0.01in}m}\hspace{0.01in}q^{-2}\psi_{L}^{\ast
}\triangleleft\partial_{A}\hspace{0.01in}(\mathcal{L}_{\partial}){^{C}%
}_{\hspace{-0.01in}B}\hspace{0.02in}\partial_{C}\circledast\psi_{R}\nonumber\\
&  -\frac{1}{4\hspace{0.01in}m}\hspace{0.01in}\psi_{L}^{\ast}\triangleleft
\partial_{A}\hspace{0.01in}(\mathcal{L}_{\partial}){^{C}}_{\hspace{-0.01in}%
B}\circledast\partial_{C}\triangleright\psi_{R}. \label{SpaTenDreDim}%
\end{align}

We can obtain the continuity equation for $i_{A}^{\ast}$ by conjugating
Eq.~(\ref{ImpDichAblSpaTen}). This way, it holds%
\begin{equation}
(i^{\ast})^{A}\,\bar{\triangleleft}\,\partial_{t}=-\hspace{0.01in}(T^{\ast
})^{AB}\,\bar{\triangleleft}\,\partial_{B}-(f^{\ast})^{A}
\label{KonKonGleImpDic}%
\end{equation}
with%
\begin{equation}
\overline{f_{A}}=(f^{\ast})^{A},\qquad\overline{T_{BA}}=(T^{\ast})^{AB}.
\end{equation}
Explicitly, we have%
\begin{equation}
(f^{\ast})^{A}=-\hspace{0.01in}\frac{1}{2}\big (\psi_{L}\circledast
(\partial^{A}\,\bar{\triangleright}\,V)\circledast\psi_{R}^{\ast}-\psi
_{L}\circledast(V\,\bar{\triangleleft}\,\partial^{A})\circledast\psi_{R}%
^{\ast}\big ) \label{KraKinIm}%
\end{equation}
and%
\begin{align}
(T^{\ast})^{AB}=  &  -\frac{1}{4\hspace{0.01in}m}\hspace{0.01in}q^{-2}%
g_{CF}\hspace{0.02in}\psi_{L}\,\bar{\triangleleft}\,\partial^{A}%
\circledast\partial^{C}(\mathcal{\bar{L}}_{\partial}){^{F}}_{\hspace
{-0.01in}E}\,\bar{\triangleright}\,\psi_{R}^{\ast}\hspace{0.02in}%
g^{EB}\nonumber\\
&  -\frac{1}{4\hspace{0.01in}m}\hspace{0.01in}g_{CF}\hspace{0.02in}\psi
_{L}\,\bar{\triangleleft}\,\partial^{A}\partial^{C}\circledast(\mathcal{\bar
{L}}_{\partial}){^{F}}_{\hspace{-0.01in}E}\,\bar{\triangleright}\,\psi
_{R}^{\ast}\hspace{0.02in}g^{EB}\nonumber\\
&  -\frac{1}{4\hspace{0.01in}m}\hspace{0.01in}q^{-2}g_{CF}\hspace{0.02in}%
\psi_{L}\circledast\partial^{C}(\mathcal{\bar{L}}_{\partial}){^{F}}%
_{\hspace{-0.01in}E}\hspace{0.02in}\partial^{A}\,\bar{\triangleright}%
\,\psi_{R}^{\ast}\hspace{0.02in}g^{EB}\nonumber\\
&  -\frac{1}{4\hspace{0.01in}m}\hspace{0.01in}g_{CF}\hspace{0.02in}\psi
_{L}\,\bar{\triangleleft}\,\partial^{C}\circledast(\mathcal{\bar{L}}%
_{\partial}){^{F}}_{\hspace{-0.01in}E}\hspace{0.02in}\partial^{A}%
\,\bar{\triangleright}\,\psi_{R}^{\ast}\hspace{0.02in}g^{EB}.
\label{StrTenImpSte}%
\end{align}

If we integrate both sides of Eq.~(\ref{ImpDichAblSpaTen}) over all space, we
obtain Newton's second law as part of the Ehrenfest theorem [also see
Eq.~(\ref{HeiGleP1Hab}) of Chap.~\ref{KapBwgGleObs}]:%
\begin{equation}
\partial_{t}\triangleright\hspace{-0.01in}\langle\hspace{0.02in}p^{A}%
\rangle=\int\text{d}_{q}^{3}\hspace{0.01in}x\,\partial_{t}\triangleright
i^{A}=\int\text{d}_{q}^{3}\hspace{0.01in}x\,f^{A}=-\langle\partial
^{A}\triangleright V\rangle.
\end{equation}
In the last step, we have used the identity%
\begin{equation}
f^{A}=\psi_{L}^{\ast}\circledast(\partial^{A}\triangleright V)\circledast
\psi_{R},
\end{equation}
which follows from Eq.~(\ref{KraDicDreDim}) by taking into acount the
following calculation [also see Eq.~(\ref{RechtsLinksDarN}) of
Chap.~\ref{KapHofStr} and Eq.~(\ref{TrivBraiV})]:%
\begin{equation}
V\triangleleft\partial^{A}=-\hspace{0.01in}\partial^{F}\triangleright
(V\triangleleft(\mathcal{L}_{\partial}){^{A}}_{\hspace{-0.01in}F})=-\delta
_{F}^{A}\,\partial^{F}\triangleright V=-\hspace{0.01in}\partial^{A}%
\triangleright V.
\end{equation}

Next, let us consider a charged particle moving in a magnetic field. To this
end, we replace the partial derivatives $\partial^{C}$ by the operators
$D^{C}$ [cf. Eq.~(\ref{DefKanImp}) of Chap.~\ref{KapHerSchGle}] and obtain the
following expressions from Eq.~(\ref{FreImpDicDefN}):%
\begin{align}
i^{C}  &  =\frac{1}{2\text{i}}\left(  \psi_{L}^{\ast}\circledast D^{C}%
\hspace{-0.01in}\triangleright\psi_{R}+\psi_{L}^{\ast}\triangleleft
\hspace{0.01in}D^{C}\hspace{-0.01in}\circledast\psi_{R}\right) \nonumber\\
&  =\frac{1}{2\text{i}}\left(  \psi_{L}^{\ast}\circledast\partial^{C}%
\hspace{-0.01in}\triangleright\psi_{R}+\psi_{L}^{\ast}\triangleleft
\hspace{0.01in}\partial^{C}\hspace{-0.01in}\circledast\psi_{R}\right)
-\psi_{L}^{\ast}\circledast e\hspace{0.01in}A^{C}\hspace{-0.01in}%
\circledast\psi_{R},\label{StrDicVekPot}\\[0.05in]
i_{C}^{\ast}  &  =\frac{1}{2\text{i}}\left(  \psi_{L}\circledast D_{C}%
\,\bar{\triangleright}\,\psi_{R}^{\ast}+\psi_{L}\,\bar{\triangleleft}%
\,\hspace{0.01in}D_{C}\circledast\psi_{R}^{\ast}\right) \nonumber\\
&  =\frac{1}{2\text{i}}\left(  \psi_{L}\circledast\partial_{C}\,\bar
{\triangleright}\,\psi_{R}^{\ast}+\psi_{L}\,\bar{\triangleleft}\,\hspace
{0.01in}\partial_{C}\circledast\psi_{R}^{\ast}\right)  -\psi_{L}\circledast
e\hspace{0.01in}A_{C}\circledast\psi_{R}^{\ast}. \label{StrDicVekPot2}%
\end{align}
Accordingly, the Hamilton operator takes on the following form:%
\begin{equation}
H=-(2\hspace{0.01in}m)^{-1}\hspace{0.01in}D^{C}D_{C}+V. \label{HamOpeVekPot}%
\end{equation}

Again, we use the Schr\"{o}\-dinger equations from Eq.~(\ref{SchGleQDef1N}) of
Chap.~\ref{KapHerSchGle} to calculate the time derivative of the momentum
density $i^{C}$ given in Eq.~(\ref{StrDicVekPot}):%
\begin{align}
\partial_{t}\triangleright i^{C}\hspace{-0.01in}=  &  \hspace{0.03in}\frac
{1}{2}\big (\psi_{L}^{\ast}\triangleleft H\circledast D^{C}\triangleright
\psi_{R}-\psi_{L}^{\ast}\circledast D^{C}\hspace{-0.01in}\triangleright
(H\triangleright\psi_{R})\big )\nonumber\\
&  \,+\frac{1}{2}\big ((\psi_{L}^{\ast}\triangleleft H)\triangleleft
D^{C}\hspace{-0.01in}\circledast\psi_{R}-\psi_{L}^{\ast}\triangleleft
\hspace{0.01in}D^{C}\hspace{-0.01in}\circledast H\triangleright\psi
_{R}\big )\nonumber\\
&  \,-\psi_{L}^{\ast}\circledast e\hspace{0.01in}\partial_{t}\triangleright
A^{C}\circledast\psi_{R}. \label{ZeiAblWskVekPot}%
\end{align}
With Eq.~(\ref{HamOpeVekPot}), we can write the first expression on the
left-hand side of the above equation as follows:%
\begin{align}
&  \psi_{L}^{\ast}\triangleleft H\circledast D^{C}\hspace{-0.01in}%
\triangleright\psi_{R}-\psi_{L}^{\ast}\circledast D^{C}\hspace{-0.01in}%
\triangleright(H\triangleright\psi_{R})=\nonumber\\
&  \qquad\qquad=-\,\psi_{L}^{\ast}\triangleleft\hspace{0.01in}D^{F}%
D_{F}\hspace{0.01in}(2\hspace{0.01in}m)^{-1}\hspace{-0.01in}\circledast
D^{C}\hspace{-0.01in}\triangleright\psi_{R}+\psi_{L}^{\ast}\circledast
(2\hspace{0.01in}m)^{-1}\hspace{0.01in}D^{C}D^{F}D_{F}\triangleright\psi
_{R}\nonumber\\
&  \qquad\qquad\hspace{0.15in}+\psi_{L}^{\ast}\circledast V\circledast
D^{C}\hspace{-0.01in}\triangleright\psi_{R}-\psi_{L}^{\ast}\circledast
D^{C}\hspace{-0.01in}\triangleright(V\circledast\psi_{R})\nonumber\\
&  \qquad\qquad=-\,\psi_{L}^{\ast}\triangleleft\hspace{0.01in}D^{F}%
D_{F}\hspace{0.01in}(2\hspace{0.01in}m)^{-1}\hspace{-0.01in}\circledast
D^{C}\hspace{-0.01in}\triangleright\psi_{R}+\psi_{L}^{\ast}\circledast
(2\hspace{0.01in}m)^{-1}\hspace{0.01in}D^{F}D_{F}D^{C}\hspace{-0.01in}%
\triangleright\psi_{R}\nonumber\\
&  \qquad\qquad\hspace{0.15in}+\psi_{L}^{\ast}\circledast f_{\text{Lor}}%
^{C}\triangleright\psi_{R}-\psi_{L}^{\ast}\circledast(\partial^{C}%
\triangleright V)\circledast\psi_{R}. \label{UmfTotDivZeiStrVek}%
\end{align}
The last expression in Eq.~(\ref{UmfTotDivZeiStrVek}) is a consequence of the
commutation relations [also see Eq.~(\ref{VerRelQuaKanImpKanImp}) of
App.~\ref{AppLorFor}]%
\begin{equation}
\lbrack D^{F}D_{F},D^{C}]=-\hspace{0.01in}2\hspace{0.01in}m\hspace
{0.01in}f_{\text{Lor}}^{C}%
\end{equation}
and%
\begin{equation}
\left[  D^{C},V\right]  =\left[  \partial^{C},V\right]  +\left[
\text{i}eA^{C},V\right]  =\partial^{C}\triangleright V.
\end{equation}
The operator for the Lorentz force density takes on the following form [cf.
Eq.~(\ref{LorKraDic}) of App.~\ref{AppLorFor}]:%
\begin{equation}
f_{\text{Lor}}^{D}=\frac{e}{2\hspace{0.01in}m\hspace{0.01in}}\text{\hspace
{0.01in}}g^{DG}\varepsilon_{ACG}\left(  \hspace{0.01in}D^{C}\hspace
{0.01in}B^{A}-B^{C}\hspace{0.01in}D^{A}\right)  .
\end{equation}
There are terms in the last expression of Eq.~(\ref{UmfTotDivZeiStrVek}) which
we can write as divergence if we apply Eq.~(\ref{GreIdeKinImp}) of
Chap.~\ref{KapGre}\ with the following identifications:%
\begin{equation}
\psi=\psi_{L}^{\ast}\hspace{0.01in}(2\hspace{0.01in}m)^{-1},\text{\qquad}%
\phi=D^{C}\hspace{-0.01in}\triangleright\psi_{R}.
\end{equation}
This way, we obtain:%
\begin{align}
&  -\psi_{L}^{\ast}\triangleleft\hspace{0.01in}D^{B}D_{B}\hspace
{0.01in}(2\hspace{0.01in}m)^{-1}\hspace{-0.01in}\circledast D^{C}%
\hspace{-0.01in}\triangleright\psi_{R}+\psi_{L}^{\ast}\circledast
(2\hspace{0.01in}m)^{-1}\hspace{0.01in}D^{B}D_{B}D^{C}\hspace{-0.01in}%
\triangleright\psi_{R}=\nonumber\\
&  \qquad\qquad=\partial^{F}\triangleright\big [q^{-2}\psi_{L}^{\ast
}\triangleleft(\mathcal{L}_{\partial}){^{B}}_{\hspace{-0.01in}F}\,D_{B}%
\hspace{0.01in}(2\hspace{0.01in}m)^{-1}\hspace{-0.01in}\circledast
D^{C}\hspace{-0.01in}\triangleright\psi_{R}\big ]\nonumber\\
&  \qquad\qquad\hspace{0.15in}+\partial^{F}\triangleright\big [\psi_{L}^{\ast
}\triangleleft(\mathcal{L}_{\partial}){^{B}}_{\hspace{-0.01in}F}%
\circledast(2\hspace{0.01in}m)^{-1}\hspace{0.01in}D_{B}D^{C}\hspace
{-0.01in}\triangleright\psi_{R}\big ]. \label{ZwiRecImpDicVek}%
\end{align}
Similar reasonings hold for the second expression on the right-hand side of
Eq.~(\ref{ZeiAblWskVekPot}):%
\begin{align}
&  (\psi_{L}^{\ast}\triangleleft H)\triangleleft\hspace{0.01in}D^{C}%
\hspace{-0.01in}\circledast\psi_{R}-\psi_{L}^{\ast}\triangleleft
\hspace{0.01in}D^{C}\circledast H\triangleright\psi_{R}=\nonumber\\
&  \qquad=\partial^{F}\triangleright\left[  q^{-2}\psi_{L}^{\ast}%
\triangleleft\hspace{0.01in}D^{C}(\mathcal{L}_{\partial}){^{B}}_{\hspace
{-0.01in}F}\,D_{B}\hspace{0.01in}(2\hspace{0.01in}m)^{-1}\hspace
{-0.01in}\circledast\psi_{R}\right] \nonumber\\
&  \qquad\hspace{0.17in}+\partial^{F}\triangleright\left[  \psi_{L}^{\ast
}\triangleleft\hspace{0.01in}D^{C}(\mathcal{L}_{\partial}){^{B}}%
_{\hspace{-0.01in}F}\hspace{0.01in}(2\hspace{0.01in}m)^{-1}\hspace
{-0.01in}\circledast D_{B}\triangleright\psi_{R}\right] \nonumber\\
&  \qquad\hspace{0.17in}-\psi_{L}^{\ast}\triangleleft f_{\text{Lor}}%
^{C}\circledast\psi_{R}+\psi_{L}^{\ast}\circledast(V\hspace{-0.01in}%
\triangleleft\partial^{C})\circledast\psi_{R}.
\end{align}
Summarizing the results so far, we finally get%
\begin{equation}
\partial_{t}\triangleright i_{C}=-\hspace{0.01in}\partial^{F}\triangleright
\hspace{0.01in}T_{FC}\hspace{0.01in}+f_{C} \label{KonGleImpDicVekPot}%
\end{equation}
with%
\begin{align}
T_{FC}=  &  -\frac{1}{4\hspace{0.01in}m}\hspace{0.01in}q^{-2}\psi_{L}^{\ast
}\triangleleft(\mathcal{L}_{\partial}){^{B}}_{\hspace{-0.01in}F}\,D_{B}%
\hspace{-0.01in}\circledast D_{C}\triangleright\psi_{R}\nonumber\\
&  -\frac{1}{4\hspace{0.01in}m}\hspace{0.01in}\psi_{L}^{\ast}\triangleleft
(\mathcal{L}_{\partial}){^{B}}_{\hspace{-0.01in}F}\circledast D_{B}%
D_{C}\triangleright\psi_{R}\nonumber\\
&  -\frac{1}{4\hspace{0.01in}m}\hspace{0.01in}q^{-2}\psi_{L}^{\ast
}\triangleleft\hspace{0.01in}D_{C}\hspace{0.01in}(\mathcal{L}_{\partial}%
){^{B}}_{\hspace{-0.01in}F}\,D_{B}\circledast\psi_{R}\nonumber\\
&  -\frac{1}{4\hspace{0.01in}m}\hspace{0.01in}\psi_{L}^{\ast}\triangleleft
\hspace{0.01in}D_{C}\hspace{0.01in}(\mathcal{L}_{\partial}){^{B}}%
_{\hspace{-0.01in}F}\circledast D_{B}\triangleright\psi_{R}
\label{SpaTenVekPot}%
\end{align}
and%
\begin{align}
f^{C}=  &  -\frac{1}{2}\big (\psi_{L}^{\ast}\circledast(\partial
^{C}\triangleright V)\circledast\psi_{R}-\psi_{L}^{\ast}\circledast
(V\hspace{-0.01in}\triangleleft\partial^{C})\circledast\psi_{R}%
\big )\nonumber\\
&  +\frac{1}{2}\big (\psi_{L}^{\ast}\circledast f_{\text{Lor}}^{C}%
\triangleright\psi_{R}-\psi_{L}^{\ast}\triangleleft f_{\text{Lor}}%
^{C}\circledast\psi_{R}\big )-\psi_{L}^{\ast}\circledast(e\hspace
{0.01in}\partial_{t}\triangleright A^{C})\circledast\psi_{R}.
\label{KraDicVekPot}%
\end{align}
Note that we can directly obtain the new stress tensor in
Eq.~(\ref{SpaTenVekPot}) from that in Eq.~(\ref{SpaTenDreDim}) by replacing
the partial derivatives $\partial^{C}$ with the operators $D^{C}$. If we
integrate the continuity equation (\ref{KonGleImpDicVekPot}) over all space
and take into account Eq.~(\ref{PatIntUneRaumInt}) of\ Chap.~\ref{KapParDer},
we get the following evolution equation [also see Eq.~(\ref{EhrTheKur1}) of
Chap.~\ref{KapBwgGleObs}]:%
\begin{equation}
\partial_{t}\triangleright\hspace{-0.01in}\langle\hspace{0.02in}p^{C}%
\rangle=-\langle\partial^{C}\triangleright V\rangle+\langle f_{\text{Lor}}%
^{C}\rangle-\langle e\hspace{0.02in}\partial_{t}\triangleright A^{C}\rangle.
\label{TimEvoImp}%
\end{equation}

By conjugating Eq.~(\ref{KonGleImpDicVekPot}), we can again obtain a second
continuity equation for momentum density [also see Eq.~(\ref{KonKonGleImpDic}%
)]. The corresponding expression for the stress tensor again follows from that
of Eq.~(\ref{StrTenImpSte}) by replacing the partial derivatives $\partial
^{C}$ with the operators $D^{C}$. For the operator of the Lorentz force
density, we now have%
\begin{equation}
(f_{\text{Lor}}^{\ast})_{D}=\frac{e}{2\hspace{0.01in}m\hspace{0.01in}%
}\text{\hspace{0.01in}}g_{DG}\hspace{0.02in}\varepsilon^{GCA}\left(
\hspace{0.01in}D_{A}\hspace{0.01in}B_{C}^{\ast}-B_{A}^{\ast}\hspace
{0.01in}D_{C}\right)
\end{equation}
with%
\begin{equation}
(B^{\ast})^{F}=\overline{B_{F}}=\text{i}\hspace{0.01in}\partial_{C}%
\,\bar{\triangleright}\,A_{D}\hspace{0.01in}\varepsilon^{FDC}\hspace
{-0.02in}=\text{i}\hspace{0.01in}A_{C}\,\bar{\triangleleft}\,\partial
_{D}\hspace{0.01in}\varepsilon^{FDC}.
\end{equation}

\section{Conservation of energy\label{ConEneKap}}

In this chapter, we derive continuity equations for the energy density of a
nonrelativistic particle. Once again, we first consider a nonrelativistic
particle in an external force field with a scalar potential $V(\mathbf{x})$.
In this case, the energy density takes on the following form:%
\begin{equation}
\mathcal{H}=\mathcal{-\hspace{0.01in}}\psi_{L}^{\ast}\hspace{-0.01in}%
\triangleleft\partial^{A}\hspace{-0.01in}\circledast(2\hspace{0.01in}%
m)^{-1}\partial_{A}\triangleright\psi_{R}+\psi_{L}^{\ast}\circledast
V\circledast\psi_{R}. \label{EneDen}%
\end{equation}

To calculate the time derivative of the energy density, we apply the
Schr\"{o}\-dinger equations given in Eq.~(\ref{SchGleQDef1N}) of
Chap.~\ref{KapHerSchGle}:%
\begin{align}
\partial_{t}\hspace{0.01in}\triangleright\mathcal{H}= &  \,\mathcal{-\hspace
{0.03in}}\text{i\hspace{0.01in}}\mathcal{(}\psi_{L}^{\ast}\triangleleft
H)\triangleleft\partial^{A}\circledast(2\hspace{0.01in}m)^{-1}\hspace
{0.01in}\partial_{A}\triangleright\psi_{R}\nonumber\\
&  +\text{i\hspace{0.01in}}\psi_{L}^{\ast}\triangleleft\hspace{0.01in}%
\partial^{A}\circledast(2\hspace{0.01in}m)^{-1}\hspace{0.01in}\partial
_{A}\triangleright(H\triangleright\psi_{R})\nonumber\\
&  +\text{i\hspace{0.01in}}\psi_{L}^{\ast}\triangleleft H\circledast
V\circledast\psi_{R}-\text{i\hspace{0.01in}}\psi_{L}^{\ast}\circledast
V\circledast H\triangleright\psi_{R}.
\end{align}
Inserting the expression for the Hamilton operator [cf. Eq.~(\ref{HamOpePot})
in Chap.~\ref{KapHerSchGle}], we obtain:%
\begin{align}
\partial_{t}\hspace{0.01in}\triangleright\mathcal{H}= &  \,\text{i\hspace
{0.01in}}\psi_{L}^{\ast}\hspace{-0.01in}\triangleleft\partial^{A}%
\partial^{\text{\hspace{0.01in}}B}\hspace{-0.01in}\partial_{B}\circledast
(2\hspace{0.01in}m)^{-2}\hspace{0.01in}\partial_{A}\hspace{0.01in}%
\triangleright\psi_{R}\nonumber\\
&  -\text{i\hspace{0.01in}}\psi_{L}^{\ast}\hspace{-0.01in}\triangleleft
\partial^{A}\circledast(2\hspace{0.01in}m)^{-2}\hspace{0.01in}\partial
^{\text{\hspace{0.01in}}B}\hspace{-0.01in}\partial_{B}\hspace{0.01in}%
\partial_{A}\triangleright\psi_{R}\nonumber\\
&  -\text{i\hspace{0.01in}}(\psi_{L}^{\ast}\circledast V)\triangleleft
\partial^{A}\circledast(2\hspace{0.01in}m)^{-1}\hspace{0.01in}\partial
_{A}\triangleright\psi_{R}\nonumber\\
&  +\text{i\hspace{0.01in}}\psi_{L}^{\ast}\hspace{-0.01in}\triangleleft
\partial^{A}\circledast(2\hspace{0.01in}m)^{-1}\hspace{0.01in}\partial
_{A}\triangleright(V\circledast\psi_{R})\nonumber\\
&  -\text{i\hspace{0.01in}}\psi_{L}^{\ast}\hspace{-0.01in}\triangleleft
\partial^{\text{\hspace{0.01in}}B}\hspace{-0.01in}\partial_{B}\hspace
{0.01in}(2\hspace{0.01in}m)^{-1}\circledast V\circledast\psi_{R}\nonumber\\
&  +\text{i\hspace{0.01in}}\psi_{L}^{\ast}\circledast V\circledast
(2\hspace{0.01in}m)^{-1}\hspace{0.01in}\partial^{\text{\hspace{0.01in}}%
B}\hspace{-0.01in}\partial_{B}\triangleright\psi_{R}.\label{ZeiAblEneDic1}%
\end{align}
We can write the first two terms on the right-hand side of the above equation
as divergence by using Eq.~(\ref{QGreIde0}) in Chap.~\ref{KapGre}:%
\begin{align}
&  \psi_{L}^{\ast}\hspace{-0.01in}\triangleleft\partial^{A}\partial
^{\text{\hspace{0.01in}}B}\hspace{-0.01in}\partial_{B}\circledast
(2\hspace{0.01in}m)^{-2}\partial_{A}\triangleright\psi_{R}-\psi_{L}^{\ast
}\hspace{-0.01in}\triangleleft\partial^{A}\circledast(2\hspace{0.01in}%
m)^{-2}\hspace{0.01in}\partial^{\text{\hspace{0.01in}}B}\hspace{-0.01in}%
\partial_{B}\hspace{0.01in}\partial_{A}\triangleright\psi_{R}=\nonumber\\
&  \qquad\qquad=-\hspace{0.01in}\partial^{C}\triangleright\big [\hspace
{0.02in}q^{-2}\hspace{0.01in}\psi_{L}^{\ast}\hspace{-0.01in}\triangleleft
\partial^{A}\hspace{0.01in}(\mathcal{L}_{\partial}){^{B}}_{\hspace{-0.01in}%
C}\text{\thinspace}\partial_{B}\circledast(2\hspace{0.01in}m)^{-2}%
\hspace{0.01in}\partial_{A}\triangleright\psi_{R}\nonumber\\
&  \qquad\qquad\hspace{0.15in}+\psi_{L}^{\ast}\hspace{-0.01in}\triangleleft
\partial^{A}(\mathcal{L}_{\partial}){^{B}}_{\hspace{-0.01in}C}\circledast
(2\hspace{0.01in}m)^{-2}\hspace{0.01in}\partial_{B}\hspace{0.01in}\partial
_{A}\hspace{0.01in}\triangleright\psi_{R}\hspace{0.02in}\big ].
\end{align}
On the right-hand side of Eq.~(\ref{ZeiAblEneDic1}), we can also combine the
third and last term as well as the fourth and penultimate term into a
divergence. However, this requires rewriting the third term with the help of
Eq.~(\ref{UmgLeiReg1}) from Chap.~\ref{KapGre} as follows:%
\begin{align}
&  -(\psi_{L}^{\ast}\circledast V)\triangleleft\partial^{A}\circledast
(2\hspace{0.01in}m)^{-1}\hspace{0.01in}\partial_{A}\triangleright\psi
_{R}=\nonumber\\
&  \qquad\qquad=-\hspace{0.01in}\psi_{L}^{\ast}\circledast V\circledast
(2\hspace{0.01in}m)^{-1}\hspace{0.01in}\partial^{A}\partial_{A}\triangleright
\psi_{R}\nonumber\\
&  \qquad\qquad\hspace{0.15in}+\partial^{C}\triangleright\hspace
{-0.01in}\big [\hspace{0.01in}(\psi_{L}^{\ast}\circledast V)\triangleleft
(\mathcal{L}_{\partial}){^{A}}_{\hspace{-0.01in}C}\circledast(2\hspace
{0.01in}m)^{-1}\hspace{0.01in}\partial_{A}\triangleright\psi_{R}%
\hspace{0.01in}\big ].
\end{align}
We do the same with the fourth term [also see Eq.~(\ref{VerAplLMat}) of
Chap.~\ref{KapGre}]:%
\begin{align}
&  \psi_{L}^{\ast}\hspace{-0.01in}\triangleleft\partial^{A}\circledast
(2\hspace{0.01in}m)^{-1}\hspace{0.01in}\partial_{A}\hspace{0.01in}%
\triangleright(V\circledast\psi_{R})=\nonumber\\
&  \qquad\qquad=\psi_{L}^{\ast}\hspace{-0.01in}\triangleleft\partial
^{A}\partial_{A}\hspace{0.01in}(2\hspace{0.01in}m)^{-1}\circledast
V\circledast\psi_{R}\nonumber\\
&  \qquad\qquad\hspace{0.15in}+\partial^{C}\triangleright\hspace
{-0.01in}\big [\hspace{0.01in}q^{-2}\psi_{L}^{\ast}\triangleleft
\hspace{0.01in}(\mathcal{L}_{\partial}){^{A}}_{\hspace{-0.01in}C}%
\,\partial_{A}\hspace{0.01in}(2\hspace{0.01in}m)^{-1}\circledast
V\circledast\psi_{R}\hspace{0.01in}\big ].
\end{align}
Inserting these results into Eq.~(\ref{ZeiAblEneDic1}), we finally get%
\begin{equation}
\partial_{t}\hspace{0.01in}\triangleright\mathcal{H}+\partial^{C}%
\triangleright S_{C}=0\label{KonEngDic}%
\end{equation}
with the following current density:%
\begin{align}
S_{C}= &  -\mathcal{\hspace{0.01in}}\frac{\text{i}}{2\hspace{0.01in}m}%
\hspace{0.01in}q^{-2}\hspace{0.01in}\psi_{L}^{\ast}\triangleleft
(\mathcal{L}_{\partial}){^{A}}_{\hspace{-0.01in}C}\hspace{0.01in}\partial
_{A}\hspace{0.01in}\circledast V\circledast\psi_{R}\nonumber\\
&  -\mathcal{\hspace{0.01in}}\frac{\text{i}}{2\hspace{0.01in}m}\hspace
{0.01in}(\psi_{L}^{\ast}\circledast V)\triangleleft(\mathcal{L}_{\partial
}){^{A}}_{\hspace{-0.01in}C}\circledast\partial_{A}\triangleright\psi
_{R}\nonumber\\
&  +\frac{\text{i}}{4\hspace{0.01in}m^{2}}\hspace{0.01in}q^{-2}\hspace
{0.01in}\psi_{L}^{\ast}\hspace{-0.01in}\triangleleft\partial^{A}%
(\mathcal{L}_{\partial}){^{B}}_{\hspace{-0.01in}C}\hspace{0.01in}\partial
_{B}\circledast\partial_{A}\triangleright\psi_{R}\nonumber\\
&  +\frac{\text{i}}{4\hspace{0.01in}m^{2}}\hspace{0.01in}\psi_{L}^{\ast
}\hspace{-0.01in}\triangleleft\partial^{A}(\mathcal{L}_{\partial}){^{B}%
}_{\hspace{-0.01in}C}\circledast\partial_{B}\hspace{0.01in}\partial
_{A}\triangleright\psi_{R}.\label{EneStrSich}%
\end{align}

If we integrate Eq.~(\ref{KonEngDic}) over all space, the divergence leads to
a surface term at spatial infinity. However, this surface term will be zero
since the wave functions vanish at spatial infinity. Accordingly, the total
energy of the nonrelativistic particle is constant over time:%
\begin{equation}
\partial_{t}\triangleright\hspace{-0.02in}\int\text{d}_{q}^{3}\hspace
{0.01in}x\,\mathcal{H}(\mathbf{x},t)=-\int\hspace{-0.01in}\text{d}_{q}%
^{3}\hspace{0.01in}x\,\partial^{C}\triangleright S_{C}=0.
\end{equation}

By conjugating Eqs.~(\ref{EneDen}) and (\ref{KonEngDic}), we get another
expression for the energy density and a corresponding continuity equation,
i.~e.%
\begin{equation}
\mathcal{H}^{\ast}=\mathcal{-\hspace{0.01in}}\psi_{L}\,\bar{\triangleleft
}\,\hspace{0.01in}\partial^{A}\hspace{0.01in}(2\hspace{0.01in}m)^{-1}%
\hspace{-0.01in}\circledast\partial_{A}\,\bar{\triangleright}\,\psi_{R}^{\ast
}+\psi_{L}\circledast V\circledast\psi_{R}^{\ast}%
\end{equation}
and%
\begin{equation}
\mathcal{H}^{\ast}\,\bar{\triangleleft}\,\partial_{t}+(S^{\ast})^{C}%
\,\bar{\triangleleft}\,\partial_{C}=0
\end{equation}
with%
\begin{align}
(S^{\ast})^{C}=  &  -\frac{\text{i}}{2\hspace{0.01in}m}\hspace{0.01in}%
q^{-2}\hspace{0.01in}\psi_{L}\circledast V\circledast g_{AB}\hspace
{0.01in}\partial^{A}\hspace{0.01in}(\mathcal{\bar{L}}_{\partial}){^{B}%
}_{\text{\hspace{-0.01in}}E}\,\bar{\triangleright}\,\psi_{R}^{\ast}%
\hspace{0.02in}g^{EC}\nonumber\\
&  -\frac{\text{i}}{2\hspace{0.01in}m}\hspace{0.01in}\psi_{L}\,\bar
{\triangleleft}\,\partial^{A}\circledast g_{AB}\hspace{0.01in}(\mathcal{\bar
{L}}_{\partial}){^{B}}_{\text{\hspace{-0.01in}}E}\,\bar{\triangleright
}\,(V\circledast\psi_{R}^{\ast})\hspace{0.02in}g^{EC}\nonumber\\
&  +\frac{\text{i}}{4\hspace{0.01in}m^{2}}\hspace{0.01in}q^{-2}\hspace
{0.01in}\psi_{L}\hspace{0.01in}\bar{\triangleleft}\,\partial^{A}\circledast
g_{BD}\hspace{0.01in}\partial^{B}(\mathcal{\bar{L}}_{\partial}){^{D}%
}_{\text{\hspace{-0.01in}}E}\hspace{0.01in}\partial_{A}\,\bar{\triangleright
}\,\psi_{R}^{\ast}\hspace{0.02in}g^{EC}\nonumber\\
&  +\frac{\text{i}}{4\hspace{0.01in}m^{2}}\hspace{0.01in}\psi_{L}%
\hspace{0.01in}\bar{\triangleleft}\,\partial^{A}\partial^{B}\circledast
g_{BD}\hspace{0.01in}(\mathcal{\bar{L}}_{\partial}){^{D}}_{\text{\hspace
{-0.01in}}E}\hspace{0.01in}\partial_{A}\,\bar{\triangleright}\,\psi_{R}^{\ast
}\hspace{0.02in}g^{EC}.
\end{align}

With some modifications, the above considerations also apply to a charged
particle moving in a magnetic field. To this end, we replace the partial
derivatives $\partial^{C}$ in Eq.~(\ref{EneDen}) with the operators $D^{C}$:%
\begin{equation}
\mathcal{H}=-\hspace{0.01in}\psi_{L}^{\ast}\hspace{-0.01in}\triangleleft
D^{C}\hspace{-0.01in}\circledast(2\hspace{0.01in}m)^{-1}D_{C}\triangleright
\psi_{R}+\psi_{L}^{\ast}\circledast V\circledast\psi_{R}.
\end{equation}
Taking the time derivative of the new energy density, we obtain:%
\begin{align}
\partial_{t}\triangleright\mathcal{H}=  &  \,\text{i\hspace{0.01in}}\psi
_{L}^{\ast}\hspace{-0.01in}\circledast e\hspace{0.01in}\partial_{t}%
\triangleright A^{C}\circledast(2\hspace{0.01in}m)^{-1}D_{C}\triangleright
\psi_{R}\nonumber\\
&  +\text{i\hspace{0.01in}}\psi_{L}^{\ast}\hspace{-0.01in}\triangleleft
D^{C}(2\hspace{0.01in}m)^{-1}\circledast e\hspace{0.01in}\partial
_{t}\triangleright A_{C}\circledast\psi_{R}\nonumber\\
&  -\text{i\hspace{0.01in}}\mathcal{(}\psi_{L}^{\ast}\triangleleft
H)\triangleleft D^{C}\hspace{-0.01in}\circledast(2\hspace{0.01in}m)^{-1}%
D_{C}\hspace{-0.01in}\triangleright\psi_{R}\nonumber\\
&  +\text{i\hspace{0.01in}}\psi_{L}^{\ast}\triangleleft\hspace{0.01in}%
D^{C}\hspace{-0.01in}\circledast(2\hspace{0.01in}m)^{-1}D_{C}\triangleright
(H\triangleright\psi_{R})\nonumber\\
&  +\text{i\hspace{0.01in}}\psi_{L}^{\ast}\triangleleft H\circledast
V\circledast\psi_{R}-\text{i\hspace{0.01in}}\psi_{L}^{\ast}\circledast
V\circledast H\triangleright\psi_{R}.
\end{align}
The first two terms on the right-hand side of the above equation form a
$q$-de\-formed power density. Once again, we can turn the other expressions
into a divergence by using Eqs.~(\ref{GreIdeKinImp}) and
(\ref{UmgLeiRegKinImp}) of Chap.~\ref{KapGre}. These calculations are very
similar to the ones that led to Eq.~(\ref{KonEngDic}). Thus we finally obtain%
\begin{align}
\partial_{t}\hspace{0.01in}\triangleright\mathcal{H}+\partial^{C}%
\triangleright S_{C}=  &  \,\frac{\text{i}e}{2\hspace{0.01in}m}\hspace
{0.01in}\psi_{L}^{\ast}\hspace{-0.01in}\circledast\partial_{t}\triangleright
A^{C}\circledast D_{C}\triangleright\psi_{R}\nonumber\\
&  +\frac{\text{i}e}{2\hspace{0.01in}m}\hspace{0.01in}\psi_{L}^{\ast}%
\hspace{-0.01in}\triangleleft D^{C}\circledast\partial_{t}\triangleright
A_{C}\circledast\psi_{R}%
\end{align}
with%
\begin{align}
S_{C}=  &  -\frac{\text{i}}{2\hspace{0.01in}m}\hspace{0.01in}q^{-2}%
\hspace{0.01in}\psi_{L}^{\ast}\triangleleft(\mathcal{L}_{\partial}){^{A}%
}_{\hspace{-0.01in}C}\hspace{0.01in}D_{A}\circledast V\circledast\psi
_{R}\nonumber\\
&  -\frac{\text{i}}{2\hspace{0.01in}m}\hspace{0.01in}(\psi_{L}^{\ast
}\circledast V)\triangleleft(\mathcal{L}_{\partial}){^{A}}_{\hspace{-0.01in}%
C}\circledast D_{A}\triangleright\psi_{R}.\nonumber\\
&  +\frac{\text{i}}{4\hspace{0.01in}m^{2}}\hspace{0.01in}q^{-2}\hspace
{0.01in}\psi_{L}^{\ast}\hspace{-0.01in}\triangleleft D^{A}(\mathcal{L}%
_{\partial}){^{B}}_{\hspace{-0.01in}C}\hspace{0.01in}D_{B}\circledast
D_{A}\triangleright\psi_{R}\nonumber\\
&  +\frac{\text{i}}{4\hspace{0.01in}m^{2}}\hspace{0.01in}\psi_{L}^{\ast
}\hspace{-0.01in}\triangleleft D^{A}(\mathcal{L}_{\partial}){^{B}}%
_{\hspace{-0.01in}C}\circledast D_{B}\hspace{0.01in}D_{A}\triangleright
\psi_{R}. \label{CurDenEneMag}%
\end{align}
Note that we can get the current density in Eq.~(\ref{CurDenEneMag}) from the
expression in Eq.~(\ref{EneStrSich}) by replacing the partial derivatives
$\partial^{C}$ with the operators $D^{C}$.

\section{Gauge transformations\label{EicTraKap}}

We will show that the continuity equation in Eq.~(\ref{KonGleImpDicVekPot}) of
Chap.~\ref{KapEneImpErh}\ is invariant under the gauge transformations%
\begin{equation}
e\hspace{0.01in}\tilde{A}^{C}\hspace{-0.01in}=e\hspace{0.01in}A^{C}%
\hspace{-0.01in}+\partial^{C}\hspace{-0.01in}\triangleright e\hspace
{0.01in}\Phi,\qquad\tilde{V}=V\hspace{-0.01in}-\partial_{0}\hspace
{0.01in}\triangleright e\hspace{0.01in}\Phi, \label{EicTraVekSkaPot}%
\end{equation}
and%
\begin{align}
\tilde{\psi}_{R}(\mathbf{x},t)  &  =\exp(\text{i}e\hspace{0.01in}%
\Phi)\circledast\psi_{R}(\mathbf{x},t),\nonumber\\
\tilde{\psi}_{L}^{\ast}(\mathbf{x},t)  &  =\psi_{L}^{\ast}(\mathbf{x}%
,t)\circledast\exp(-\text{i}e\hspace{0.01in}\Phi).
\end{align}
Note that $\Phi(\mathbf{x})$ is a central element of the algebra of position
space [cf. Eq.~(\ref{CenEle}) of Chap.~\ref{KapHerSchGle}]. In addition to
this, we require that $\Phi(\mathbf{x})$ has trivial braiding.

First, we show that the momentum density in Eq.~(\ref{StrDicVekPot}) of
Chap.~\ref{KapEneImpErh} is invariant under the above gauge transformations.
To this end, we do the following calculation:%
\begin{align}
\text{i}^{-1}\tilde{D}^{C}\triangleright\tilde{\psi}_{R}=  &  \hspace
{0.04in}\text{i}^{-1}\partial^{C}\triangleright\big (\exp(\text{i}%
e\hspace{0.01in}\Phi)\circledast\psi_{R}\big )-e\hspace{0.01in}A^{C}%
\hspace{-0.01in}\circledast\exp(\text{i}e\hspace{0.01in}\Phi)\circledast
\psi_{R}\nonumber\\
&  -(\partial^{C}\hspace{-0.01in}\triangleright e\hspace{0.01in}%
\Phi)\circledast\exp(\text{i}e\hspace{0.01in}\Phi)\circledast\psi
_{R}\nonumber\\
=  &  \hspace{0.04in}\big (\text{i}^{-1}\partial^{C}\hspace{-0.01in}%
\triangleright\exp(\text{i}e\hspace{0.01in}\Phi)\big )\circledast\psi_{R}%
+\exp(\text{i}e\hspace{0.01in}\Phi)\circledast\text{i}^{-1}\partial^{C}%
\hspace{-0.01in}\triangleright\psi_{R}\nonumber\\
&  -e\hspace{0.01in}A^{C}\hspace{-0.01in}\circledast\exp(\text{i}%
e\hspace{0.01in}\Phi)\circledast\psi_{R}-(\partial^{C}\hspace{-0.01in}%
\triangleright e\hspace{0.01in}\Phi)\circledast\exp(\text{i}e\hspace
{0.01in}\Phi)\circledast\psi_{R}. \label{WirKanImpEichInv}%
\end{align}
The last step follows from the Leibniz rules of the $q$-de\-formed partial
derivatives and the trivial braiding of $\Phi(\mathbf{x})$. For the same
reasons, we also have%
\begin{align}
\partial^{C}\triangleright(e\hspace{0.01in}\Phi)^{n}  &  =\sum_{j\hspace
{0.01in}=\hspace{0.01in}0}^{n\hspace{0.01in}-1}\hspace{0.01in}(e\hspace
{0.01in}\Phi)^{j}\circledast(\partial^{C}\hspace{-0.01in}\triangleright
e\hspace{0.01in}\Phi)\circledast(e\hspace{0.01in}\Phi)^{n\hspace{0.01in}%
-1-j}\nonumber\\
&  =n\hspace{0.01in}(\partial^{C}\hspace{-0.01in}\triangleright e\hspace
{0.01in}\Phi)\circledast(e\hspace{0.01in}\Phi)^{n\hspace{0.01in}-1}
\label{ParAblSkaPotN}%
\end{align}
with%
\begin{equation}
(e\hspace{0.01in}\Phi)^{n}=\hspace{0.01in}\underset{n\text{-mal}}%
{\underbrace{e\hspace{0.01in}\Phi\circledast\ldots\circledast e\hspace
{0.01in}\Phi}}.
\end{equation}
From Eq.~(\ref{ParAblSkaPotN}) follows:%
\begin{align}
\text{i}^{-1}\partial^{C}\hspace{-0.01in}\triangleright\exp(\text{i}%
e\hspace{0.01in}\Phi)  &  =\sum_{n\hspace{0.01in}=\hspace{0.01in}0}^{\infty
}\frac{\text{i}^{n}}{n!}\,\text{i}^{-1}\partial^{C}\hspace{-0.01in}%
\triangleright(e\hspace{0.01in}\Phi)^{n}\nonumber\\
&  =\sum_{n\hspace{0.01in}=1}^{\infty}\frac{\text{i}^{n\hspace{0.01in}-1}%
}{(n-1)!}\,(\partial^{C}\hspace{-0.01in}\triangleright e\hspace{0.01in}%
\Phi)\circledast(e\hspace{0.01in}\Phi)^{n\hspace{0.01in}-1}\nonumber\\
&  =(\partial^{C}\hspace{-0.01in}\triangleright e\hspace{0.01in}%
\Phi)\circledast\sum_{n\hspace{0.01in}=\hspace{0.01in}0}^{\infty}%
\frac{\text{i}^{n}}{n!}\,(e\hspace{0.01in}\Phi)^{n}=(\partial^{C}%
\hspace{-0.01in}\triangleright e\hspace{0.01in}\Phi)\circledast\exp
(\text{i}e\hspace{0.01in}\Phi). \label{AblPhaFak}%
\end{align}
If we insert this result into Eq.~(\ref{WirKanImpEichInv}), we obtain:%
\begin{align}
\text{i}^{-1}\tilde{D}^{C}\hspace{-0.01in}\triangleright\tilde{\psi}_{R}=  &
\hspace{0.03in}\big (\text{i}^{-1}\partial^{C}\hspace{-0.01in}\triangleright
\exp(\text{i}e\hspace{0.01in}\Phi)\big )\circledast\psi_{R}+\exp
(\text{i}e\hspace{0.01in}\Phi)\circledast\text{i}^{-1}\partial^{C}%
\hspace{-0.01in}\triangleright\psi_{R}\nonumber\\
&  -e\hspace{0.01in}A^{C}\hspace{-0.01in}\circledast\exp(\text{i}%
e\hspace{0.01in}\Phi)\circledast\psi_{R}-(\partial^{C}\hspace{-0.01in}%
\triangleright e\hspace{0.01in}\Phi)\circledast\exp(\text{i}e\hspace
{0.01in}\Phi)\circledast\psi_{R}\nonumber\\
=  &  \hspace{0.03in}\exp(\text{i}e\hspace{0.01in}\Phi)\circledast
\text{i}^{-1}\partial^{C}\hspace{-0.01in}\triangleright\psi_{R}-\exp
(\text{i}e\hspace{0.01in}\Phi)\circledast e\hspace{0.01in}A^{C}\hspace
{-0.01in}\circledast\psi_{R}\nonumber\\
=  &  \hspace{0.03in}\exp(\text{i}e\hspace{0.01in}\Phi)\circledast
\text{i}^{-1}D^{C}\hspace{-0.01in}\triangleright\psi_{R}.
\label{KanImpLinWirEic}%
\end{align}
In the same manner, we get the following identity:%
\begin{equation}
\tilde{\psi}_{L}^{\ast}\triangleleft\hspace{0.01in}\tilde{D}^{C}%
\hspace{-0.01in}=\psi_{L}^{\ast}\triangleleft\hspace{0.01in}D^{C}%
\hspace{-0.01in}\circledast\exp(-\text{i}e\hspace{0.01in}\Phi).
\label{KanImpRecWirEic}%
\end{equation}
Applying Eqs.~(\ref{KanImpLinWirEic}) and (\ref{KanImpRecWirEic}), we find
that the momentum density in Eq.~(\ref{StrDicVekPot}) of
Chap.~\ref{KapEneImpErh} is invariant under gauge transformations:%
\begin{align}
\widetilde{i}^{\hspace{0.01in}C}=  &  \hspace{0.03in}\frac{1}{2\text{i}%
}\big (\tilde{\psi}_{L}^{\ast}\circledast\hspace{0.01in}\tilde{D}^{C}%
\hspace{-0.01in}\triangleright\tilde{\psi}_{R}+\tilde{\psi}_{L}^{\ast
}\triangleleft\hspace{0.01in}\tilde{D}^{C}\hspace{-0.01in}\circledast
\tilde{\psi}_{R}\big )\nonumber\\
=  &  \hspace{0.03in}\frac{1}{2\text{i}}\psi_{L}^{\ast}\circledast
\exp(-\text{i}e\hspace{0.01in}\Phi)\circledast\exp(\text{i}e\hspace
{0.01in}\Phi)\circledast D^{C}\hspace{-0.01in}\triangleright\psi
_{R}\nonumber\\
&  \hspace{0.03in}+\frac{1}{2\text{i}}\psi_{L}^{\ast}\hspace{-0.01in}%
\triangleleft D^{C}\hspace{-0.01in}\circledast\exp(-\text{i}e\hspace
{0.01in}\Phi)\circledast\exp(\text{i}e\hspace{0.01in}\Phi)\circledast\psi
_{R}\nonumber\\
=  &  \hspace{0.03in}\frac{1}{2\text{i}}\big (\psi_{L}^{\ast}\circledast
D^{C}\hspace{-0.01in}\triangleright\psi_{R}+\psi_{L}^{\ast}\hspace
{-0.01in}\triangleleft D^{C}\hspace{-0.01in}\circledast\psi_{R}\big )=i^{C}.
\end{align}

Similar considerations show that the $q$-de\-formed stress tensor in
Eq.~(\ref{SpaTenVekPot}) of Chap.~\ref{KapEneImpErh} is also invariant under
gauge transformations. That the $q$-de\-formed force density in
Eq.~(\ref{KraDicVekPot}) of Chap.~\ref{KapEneImpErh} is gauge invariant can be
proven in complete analogy to the undeformed case. This results arises from
the fact that the $q$-de\-formed magnetic field does not change under the
gauge transformations in Eq.~(\ref{EicTraVekSkaPot}):%
\begin{equation}
\tilde{B}_{F}=\text{i}\hspace{0.01in}\partial^{C}\triangleright\tilde{A}%
^{D}\varepsilon_{DCF}=\text{i}\hspace{0.01in}\partial^{C}\triangleright
A^{D}\varepsilon_{DCF}-\text{i}\hspace{0.01in}\partial^{C}\partial
^{D}\triangleright\Phi\hspace{0.02in}\varepsilon_{DCF}=B_{F}.
\end{equation}
Note that the last step in the above calculation is a consequence of the
following identity \cite{Lorek:1995ph}:%
\[
\partial^{C}\partial^{D}\varepsilon_{DCF}=0.
\]

\section{Equations of motion\label{KapBwgGleObs}}

In Ref.~\cite{Wachter:2020A}, we have shown that the Heisenberg equation of
motion also holds for quantum-mechanical systems of $q$-de\-formed Euclidean
space. Using Heisenberg's equation, we can calculate the time derivative of an
observable of a $q$-de\-formed nonrelativistic particle. We will see that
these results are in agreement with the continuity equations in
Chap.~\ref{KapEneImpErh}.

For the momentum operator of a free $q$-de\-formed nonrelativistic particle,
we get the following equation of motion in the Heisenberg picture:%
\begin{align}
\partial_{t}\hspace{0.01in}\triangleright(P^{A})_{H} &  =\text{i\hspace
{0.01in}}[H_{0},(P^{A})_{H}]\nonumber\\
&  =-\hspace{0.01in}\text{i\hspace{0.01in}}(2\hspace{0.01in}m)^{-1}%
\hspace{0.01in}\partial^{C}\partial_{C}\hspace{0.01in}\text{i}^{-1}%
\partial^{A}+\text{i}^{-1}\partial^{A}\,\text{i\hspace{0.01in}}(2\hspace
{0.01in}m)^{-1}\hspace{0.01in}\partial^{C}\partial_{C}=0.\label{HeiGleMomHab}%
\end{align}
In the above calculation, we have used the fact that the free Hamilton
operator $H_{0}$ commutes with the $q$-de\-formed partial derivatives [also
see Eq.~(\ref{ComHP}) in Chap.~\ref{KapHerSchGle}]:%
\begin{align}
(P_{A})_{H} &  =\exp(\text{i}tH_{0})\,P^{A}\hspace{0.01in}\exp(-\text{i}%
tH_{0})\nonumber\\
&  =\exp(\text{i}tH_{0})\,\text{i}^{-1}\partial^{A}\hspace{0.01in}%
\exp(-\text{i}tH_{0})=\,\text{i}^{-1}\partial^{A}=P^{A}.
\end{align}
By Eq.~(\ref{HeiGleMomHab}), we learn that the $q$-de\-formed momentum of a
free particle is time-in\-de\-pen\-dent.

Next, we calculate the time derivative of the position operator in the
Heisenberg picture:%
\begin{align}
\partial_{t}\hspace{0.01in}\triangleright(X^{A})_{H} &  =\text{i\hspace
{0.01in}}[H_{0},(X^{A})_{H}]\nonumber\\
&  =-\,\text{i\hspace{0.01in}}(2\hspace{0.01in}m)^{-1}\partial^{C}\partial
_{C}\text{\hspace{0.01in}}(X^{A})_{H}+\text{i\hspace{0.01in}}(X^{A}%
)_{H}\hspace{0.02in}(2\hspace{0.01in}m)^{-1}\partial^{C}\partial
_{C}.\label{ZwiZeiAblOrtOpeHei}%
\end{align}
To\ evaluate the commutator with the free Hamilton operator,\ we have used the
following identity\footnote{If we use the partial derivatives $\hat{\partial
}^{A}$, $\ $we must replace $q$ by $q^{-1}$ in the following formulas.}%
\begin{equation}
g_{CD}\hspace{0.01in}\partial^{C}\partial^{D}X^{A}=(1+q^{-2})\hspace
{0.01in}\partial^{A}\hspace{-0.01in}+q^{4}X^{A}g_{CD}\hspace{0.01in}%
\partial^{C}\partial^{D}.\label{VerImpQuaOrtKoo}%
\end{equation}
To get the latter, we have applied the Leibniz rules in
Eq.~(\ref{DifKalExtEukQuaDreUnk}) of Chap.~\ref{KapParDer}.

The time derivative of the position operator should be proportional to the
momentum operator. This assumption requires that mass $m$ is a scaling
operator:%
\begin{equation}
m\hspace{0.01in}X^{A}=q^{4}X^{A}m.\label{SkaMasKoo}%
\end{equation}
From Eqs.~(\ref{VerImpQuaOrtKoo}) and (\ref{SkaMasKoo}) follows:%
\begin{gather}
-\,\text{i\hspace{0.01in}}(2\hspace{0.01in}m)^{-1}\partial^{C}\partial
_{C}\text{\hspace{0.01in}}(X^{A})_{H}=\nonumber\\
=-\text{\hspace{0.01in}i\hspace{0.01in}}(1+q^{-2})(2\hspace{0.01in}%
m)^{-1}\partial^{A}\hspace{-0.01in}-\text{i\hspace{0.01in}}(X^{A}%
)_{H}\text{\hspace{0.01in}}(2\hspace{0.01in}m)^{-1}\hspace{0.01in}\partial
^{C}\partial_{C}.
\end{gather}
Inserting this result into Eq.~(\ref{ZwiZeiAblOrtOpeHei}), we finally get:%
\begin{equation}
\partial_{t}\hspace{0.01in}\triangleright(X^{A})_{H}=\text{i\hspace{0.01in}%
}[H_{0},(X^{A})_{H}]=\frac{[[2]]_{q^{-2}}}{2\hspace{0.01in}m}\hspace
{0.01in}P^{A}.\label{HeiGlImp1Hab}%
\end{equation}

Next, we will study how a time-in\-de\-pen\-dent scalar potential
$V(\mathbf{x})$ changes the Heisenberg equation of motion for the momentum
operator or the position operator. For this reason, we consider the following
Hamilton operator:%
\begin{equation}
H=H_{0}+V(\mathbf{x}).\label{HamOpeSkaBew}%
\end{equation}
The Heisenberg equation of motion for the momentum operator then reads:%
\begin{align}
\partial_{t}\hspace{0.01in}\triangleright(P^{A})_{H} &  =\text{i\hspace
{0.01in}}[H_{H},(P^{A})_{H}]=\text{i\hspace{0.01in}}[(H_{0})_{H},(P^{A}%
)_{H}]+\text{i\hspace{0.01in}}[V_{H},(P^{A})_{H}]\nonumber\\
&  =(V(\mathbf{x})\,\partial^{A})_{H}-(\partial^{A}\hspace{0.01in}%
V(\mathbf{x}))_{H}\nonumber\\
&  =(V(\mathbf{x})\,\partial^{A})_{H}-\big ([(\partial^{A})_{(1)}%
\triangleright V(\mathbf{x})]\hspace{0.01in}(\partial^{A})_{(2)}%
\big )_{H}\nonumber\\
&  =(V(\mathbf{x})\,\partial^{A})_{H}-(\partial^{A}\triangleright
V(\mathbf{x}))_{H}-(V(\mathbf{x})\,\partial^{A})_{H}\nonumber\\
&  =-\hspace{0.01in}(\partial^{A}\triangleright V(\mathbf{x}))_{H}%
.\label{HeiGleP1Hab}%
\end{align}
The above calculation employs the fact that $H_{0}$ commutes with the momentum
operator [also see Eq.~(\ref{HeiGleMomHab})] and that $V(\mathbf{x})$ has
trivial braiding properties. Note that Eq.~(\ref{HeiGleP1Hab}) yields a
$q$\textit{-ana\-log of Newton's second law}.

The Heisenberg equation of motion for the position operator remains unchanged
compared to Eq.~(\ref{HeiGlImp1Hab}) since $V(\mathbf{x})$ commutes with the
position operator:%
\begin{align}
\partial_{t}\hspace{0.01in}\triangleright(X^{A})_{H}  &  =\text{i\hspace
{0.01in}}[H_{H},(X^{A})_{H}]=\text{i\hspace{0.01in}}[(H_{0})_{H},(X^{A}%
)_{H}]+\text{i\hspace{0.01in}}[V_{H},(X^{A})_{H}]\nonumber\\
&  =\frac{[[2]]_{q^{-2}}}{2\hspace{0.01in}m}\hspace{0.01in}(P^{A})_{H}.
\label{HeiGleX1}%
\end{align}
Combining Eq.~(\ref{HeiGleP1Hab}) with Eq.~(\ref{HeiGleX1}) gives:%
\begin{equation}
(\partial_{t})^{2}\triangleright(X^{A})_{H}=-\hspace{0.01in}(\partial
^{A}\triangleright V)_{H}.
\end{equation}

We can also write down $q$\textit{-de\-formed versions of Hamilton's
equations}. To this end, we do the calculations%
\begin{equation}
\partial_{x}^{A}\triangleright H=\partial_{x}^{A}\triangleright H_{0}%
+\partial_{x}^{A}\triangleright V(\mathbf{x})=\partial_{x}^{A}\triangleright V
\end{equation}
and\footnote{The partial derivatives $\partial_{p}^{A}$ act on momentum space
in the same way as the partial derivatives $\partial_{x}^{A}$ act on position
space.}%
\begin{align}
\partial_{p}^{A}\triangleright H  &  =\partial_{p}^{A}\triangleright
H_{0}+\partial_{p}^{A}\triangleright\hspace{-0.01in}V\nonumber\\
&  =\hspace{0.01in}\partial_{p}^{A}\triangleright g_{AB}\hspace{0.01in}%
P^{A}P^{B}\text{\hspace{0.01in}}(2\hspace{0.01in}m)^{-1}=[[2]]_{q^{-2}}%
\hspace{0.01in}P^{A}\text{\hspace{0.01in}}(2\hspace{0.01in}m)^{-1}.
\end{align}
Comparing these results with those of Eqs.~(\ref{HeiGleP1Hab}) and
(\ref{HeiGleX1}), we get the following identities:%
\begin{equation}
\partial_{t}\hspace{0.01in}\triangleright(P^{A})_{H}=-\hspace{0.01in}%
(\partial_{x}^{A}\triangleright H)_{H},\qquad\partial_{t}\hspace
{0.01in}\triangleright(X^{A})_{H}=(\partial_{p}^{A}\triangleright H)_{H}.
\end{equation}

We again consider the case where the Hamiltonian operator depends not only on
a scalar potential $V(\mathbf{x})$ but also on a vector potential
$\mathbf{A}(\mathbf{x})$ [cf. Eq.~(\ref{HamVekPot}) of
Chap.~\ref{KapHerSchGle}], i.~e.%
\begin{equation}
H=(2\hspace{0.01in}m)^{-1}\text{\hspace{0.01in}}\Pi^{C}\Pi_{C}+V(\mathbf{x})
\end{equation}
with%
\begin{equation}
\Pi^{C}=P^{C}\hspace{-0.01in}-e\hspace{0.01in}A^{C}(\mathbf{x}%
).\label{kinMomDef}%
\end{equation}
We recall that the components $e\hspace{0.01in}A^{C}$ have the same braiding
properties as the components $P^{C}$ of the momentum operator. Consequently,
if the momentum operator is represented by the derivatives $\partial^{C}$, the
components $e\hspace{0.01in}A^{C}$ and the coordinate generators $X^{D}$
satisfy the following commutation relations [also see
Eq.~(\ref{DifKalExtEukQuaDreUnk})\ of Chap.~\ref{KapParDer}]:\footnote{If the
momentum operator is represented by the derivatives $\hat{\partial}^{C}$, we
must replace $\hat{R}^{-1}$ by $\hat{R}$ and $q$ by $q^{-1}$ in the following
identities [also see Eqs.~(\ref{DifKalExtEukQuaDreUnk}) and
(\ref{DifKalExtEukQuaDreKon})\ of Chap.~\ref{KapParDer}].}%
\begin{equation}
e\hspace{0.01in}A^{C}\hspace{0.01in}X^{D}=(\hat{R}^{-1}){^{CD}}_{\hspace
{-0.01in}EF}\hspace{0.01in}X^{E}e\hspace{0.01in}A^{F}.
\end{equation}
The same considerations that led to the identities in
Eq.~(\ref{VerImpQuaOrtKoo}) yield:%
\begin{align}
g_{CD}\hspace{0.01in}e\hspace{0.01in}A^{C}e\hspace{0.01in}A^{D}X^{E} &
=q^{4}X^{E}\hspace{0.01in}g_{CD}\hspace{0.01in}e\hspace{0.01in}A^{C}%
e\hspace{0.01in}A^{D}\nonumber\\
g_{CD}\hspace{0.01in}e\hspace{0.01in}A^{C}P^{D}X^{E} &  =-\text{\hspace
{0.01in}i\hspace{0.01in}}e\hspace{0.01in}A^{E}\hspace{-0.01in}+q^{4}%
X^{E}\hspace{0.01in}g_{CD}\hspace{0.01in}e\hspace{0.01in}A^{C}P^{D}%
,\nonumber\\
g_{CD}\hspace{0.01in}P^{C}e\hspace{0.01in}A^{D}X^{E} &  =-\text{\hspace
{0.01in}i\hspace{0.01in}}q^{-2}e\hspace{0.01in}A^{E}\hspace{-0.01in}%
+q^{4}X^{E}\hspace{0.01in}g_{CD}\hspace{0.01in}P^{C}e\hspace{0.01in}A^{D}.
\end{align}
Together with Eq.~(\ref{SkaMasKoo}), the identities above imply:%
\begin{align}
\lbrack(2\hspace{0.01in}m)^{-1}e\hspace{0.01in}A^{D}e\hspace{0.01in}%
A_{D},X^{E}] &  =0,\nonumber\\
\lbrack(2\hspace{0.01in}m)^{-1}e\hspace{0.01in}A^{D}P_{D},X^{E}] &
=-\text{\hspace{0.01in}i\hspace{0.01in}}(2\hspace{0.01in}m)^{-1}%
e\hspace{0.01in}A^{E},\nonumber\\
\lbrack(2\hspace{0.01in}m)^{-1}P^{D}e\hspace{0.01in}A_{D},X^{E}] &
=-\text{\hspace{0.01in}i\hspace{0.01in}}(2\hspace{0.01in}m)^{-1}q^{-2}%
e\hspace{0.01in}A^{E}.
\end{align}
We can calculate the time derivative of the position operator in the
Heisenberg picture by using the commutators above:%
\begin{align}
\partial_{t}\hspace{0.01in}\triangleright(X^{D})_{H} &  =\text{i\hspace
{0.01in}}[H_{H},(X^{D})_{H}]=[[2]]_{q^{-2}}(2\hspace{0.01in}m)^{-1}%
(P^{D}\hspace{-0.01in}-e\hspace{0.01in}A^{D})_{H}\nonumber\\
&  =\frac{[[2]]_{q^{-2}}}{2\hspace{0.01in}m}\hspace{0.01in}(\Pi^{D}%
)_{H}.\label{ZeiAblOrtVekPot}%
\end{align}

Next, we calculate the time derivative of kinetic momentum. Using
Eq.~(\ref{VerRelQuaKanImpKanImp}) in App.~\ref{AppLorFor} and the identities%
\begin{equation}
\lbrack V,\Pi^{D}]=[V,P^{D}\hspace{-0.01in}-e\hspace{0.01in}A^{D}%
]=-\hspace{0.01in}\partial^{D}\triangleright V,\qquad\frac{\partial
\hspace{0.01in}\Pi^{D}}{\partial\text{\hspace{0.01in}}t}=-e\hspace
{0.01in}\frac{\partial A^{D}}{\partial\text{\hspace{0.01in}}t},
\end{equation}
we can do the following calculation:%
\begin{align}
\partial_{t}\hspace{0.01in}\triangleright(\Pi^{D})_{H} &  =\text{i\hspace
{0.01in}}[H,(\Pi^{D})_{H}]+\frac{\partial\hspace{0.01in}(\Pi^{D})_{H}%
}{\partial\text{\hspace{0.01in}}t}\nonumber\\
&  =\text{i\hspace{0.01in}}[(2\hspace{0.01in}m)^{-1}(\Pi^{C})_{H}(\Pi_{C}%
)_{H},(\Pi^{D})_{H}]+\text{i\hspace{0.01in}}[V_{H},(\Pi^{D})_{H}%
]+\frac{\partial\hspace{0.01in}(\Pi^{D})_{H}}{\partial t}\nonumber\\
&  =(f_{\text{Lor}}^{D})_{H}-e\hspace{0.01in}\frac{\partial\hspace
{0.01in}(A^{D})_{H}}{\partial\text{\hspace{0.01in}}t}-(\partial^{D}%
\triangleright V)_{H}.\label{EhrTheKur1}%
\end{align}
Combining Eq.~(\ref{ZeiAblOrtVekPot}) with the result above, we get:%
\begin{equation}
\frac{2\hspace{0.01in}m}{[[2]]_{q^{-2}}}\hspace{0.01in}(\partial_{t}%
)^{2}\triangleright(X^{D})_{H}=(f_{\text{Lor}}^{D})_{H}-e\hspace
{0.01in}(\partial_{t}\triangleright A^{D})_{H}-(\partial^{D}\triangleright
V)_{H}.\label{NewLorKra}%
\end{equation}
By taking the expectation values of both sides, we finally obtain%
\begin{align}
\frac{2\hspace{0.01in}m}{[[2]]_{q^{-2}}}\hspace{0.01in}(\partial_{t})^{2} &
\triangleright\int\text{d}_{q}^{3}\hspace{0.01in}x\,\psi_{L}^{\ast}%
(\mathbf{x},0)\circledast(X^{D})_{H}\triangleright\psi_{R}(\mathbf{x}%
,0)=\nonumber\\
= &  \,\int\text{d}_{q}^{3}\hspace{0.01in}x\,\psi_{L}^{\ast}(\mathbf{x}%
,0)\circledast(f_{\text{Lor}}^{D})_{H}\circledast\psi_{R}(\mathbf{x}%
,0)\nonumber\\
&  -\int\text{d}_{q}^{3}\hspace{0.01in}x\,\psi_{L}^{\ast}(\mathbf{x}%
,0)\circledast e\hspace{0.01in}(\partial_{t}\triangleright A^{D}%
)_{H}\circledast\psi_{R}(\mathbf{x},0)\nonumber\\
&  -\int\text{d}_{q}^{3}\hspace{0.01in}x\,\psi_{L}^{\ast}(\mathbf{x}%
,0)\circledast(\partial^{D}\triangleright V)_{H}\circledast\psi_{R}%
(\mathbf{x},0),
\end{align}
or in abbreviated form:%
\begin{equation}
\frac{2\hspace{0.01in}m}{[[2]]_{q^{-2}}}\hspace{0.01in}(\partial_{t}%
)^{2}\triangleright\left\langle X^{D}\right\rangle _{\psi}=\langle
f_{\text{Lor}}^{D}\rangle_{\psi}-e\hspace{0.01in}\langle\partial
_{t}\triangleright A^{D}\rangle_{\psi}-\langle\partial^{D}\triangleright
V\rangle_{\psi}.
\end{equation}

Finally, we consider the time derivative of the $q$-de\-formed angular
momentum operator. We recall that $H$ is a scalar as to the actions of the
Hopf algebra $\mathcal{U}_{q}(\operatorname*{su}_{2})$. For this reason, $H$
commutes with all $g\in\mathcal{U}_{q}(\operatorname*{su}_{2})$:%
\begin{equation}
g\hspace{0.01in}H=(\hspace{0.01in}g_{(1)}\triangleright H)\hspace
{0.01in}g_{(2)}=\varepsilon(\hspace{0.01in}g_{(1)})\hspace{0.01in}%
H\hspace{0.01in}g_{(2)}=H\hspace{0.01in}g.
\end{equation}
The components of $q$-de\-formed angular momentum are elements of
$\mathcal{U}_{q}(\operatorname*{su}_{2})$ (also see App.~\ref{AppAngMom}). Due
to this fact, $H$ commutes with the $q$-de\-formed angular momentum operator:%
\begin{align}
\partial_{t}\triangleright(L^{A})_{H} &  =\text{i\hspace{0.01in}}%
[H,(L^{A})_{H}]=\text{i\hspace{0.01in}}HL^{A}\hspace{-0.01in}-\text{i\hspace
{0.01in}}L^{A}H\nonumber\\
&  =\text{i\hspace{0.01in}}HL^{A}\hspace{-0.01in}-\text{i\hspace{0.01in}%
}HL^{A}=0.\label{BewGleDreImpDre}%
\end{align}
The above result implies that the expectation value of each component of
$q$-de\-formed angular momentum is constant over time:%
\begin{align}
\partial_{t}\triangleright\langle L^{D}\rangle_{\psi} &  =\partial
_{t}\triangleright\hspace{-0.02in}\int\text{d}_{q}^{3}\hspace{0.01in}%
x\,\psi_{L}^{\ast}(\mathbf{x},0)\circledast(L^{D})_{H}\triangleright\psi
_{R}(\mathbf{x},0)\nonumber\\
&  =\int\text{d}_{q}^{3}\hspace{0.01in}x\,\psi_{L}^{\ast}(\mathbf{x}%
,0)\circledast(\partial_{t}\triangleright(L^{D})_{H})\triangleright\psi
_{R}(\mathbf{x},0)=0.
\end{align}

\appendix

\section{$q$-Analog of Lorentz force density\label{AppLorFor}}

The commutation relations between the components of $q$-de\-formed kinetic
momentum depend on the epsilon tensor of three-di\-men\-sion\-al
$q$-de\-formed Euclidean space. Its non-van\-ish\-ing components are
\cite{2004JPhA...37.9175F,Lorek:1997eh,majid-1995-36,Meyer:1994wi,Schmidt:2005ph}%
:%
\begin{align}
\varepsilon^{-3+} &  =\varepsilon_{+3-}=-q^{-2}, & \varepsilon^{3-+} &
=\varepsilon_{+-3}=1,\nonumber\\
\varepsilon^{-+3} &  =\varepsilon_{3+-}=1, & \varepsilon^{+-3} &
=\varepsilon_{3-+}=-1,\nonumber\\
\varepsilon^{3+-} &  =\varepsilon_{-+3}=-1, & \varepsilon^{+3-} &
=\varepsilon_{-3+}=q^{2},\nonumber\\
\varepsilon^{333} &  =\varepsilon_{333}=-\lambda. &  & \label{epseu3app}%
\end{align}
We can get the commutation relations for the components of \textit{kinetic
momentum} [see Eq.~(\ref{kinMomDef}) in Chap.~\ref{KapBwgGleObs}] by the
following calculation:%
\begin{align}
&  \Pi^{C}\Pi^{D}\varepsilon_{DCF}=\,(P^{C}\hspace{-0.01in}-e\hspace
{0.01in}A^{C})(P^{D}\hspace{-0.01in}-e\hspace{0.01in}A^{D})\text{\thinspace
}\varepsilon_{DCF}\nonumber\\
&  \qquad=\,P^{C}P^{D}\varepsilon_{DCF}+e\hspace{0.01in}A^{C}e\hspace
{0.01in}A^{D}\varepsilon_{DCF}-P^{C}e\hspace{0.01in}A^{D}\varepsilon
_{DCF}-e\hspace{0.01in}A^{C}P^{D}\varepsilon_{DCF}\nonumber\\
&  \qquad=-\hspace{0.01in}\text{i}^{-1}\partial^{C}\triangleright
e\hspace{0.01in}A^{D}\varepsilon_{DCF}-q^{-4}\varepsilon_{DCF}\hspace
{0.01in}(\hat{R}^{-1}){^{CD}}_{\hspace{-0.01in}GH}\,e\hspace{0.01in}A^{G}%
P^{H}-e\hspace{0.01in}A^{C}P^{D}\varepsilon_{DCF}\nonumber\\
&  \qquad=e\hspace{0.01in}B_{F}+e\hspace{0.01in}A^{G}P^{H}\varepsilon
_{HGF}-e\hspace{0.01in}A^{C}P^{D}\varepsilon_{DCF}=e\hspace{0.01in}%
B_{F}.\label{VerRelKinImp}%
\end{align}
In the third step, we used the following identities \cite{Lorek:1997eh}:%
\begin{equation}
P^{C}P^{D}\varepsilon_{DCF}=A^{C}A^{D}\varepsilon_{DCF}=0.
\end{equation}
Furthermore, we took into account that $q$-de\-formed partial derivatives give
a representation of canonical momentum:%
\begin{equation}
P^{C}e\hspace{0.01in}A^{D}=\text{i}^{-1}\partial^{C}\triangleright
e\hspace{0.01in}A^{D}+q^{-4}(\hat{R}^{-1}){^{CD}}_{\hspace{-0.01in}%
GH}\,e\hspace{0.01in}A^{G}P^{H}.\label{VerImpVekMax}%
\end{equation}
The penultimate step in Eq.~(\ref{VerRelKinImp}) follows with (also cf.
Ref.~\cite{Lorek:1997eh})%
\begin{equation}
\varepsilon_{DCF}\hspace{0.01in}(\hat{R}^{-1}){^{CD}}_{\hspace{-0.01in}%
GH}=-\hspace{0.01in}q^{4}\hspace{0.01in}\varepsilon_{HGF}%
\end{equation}
and the expression for $q$\textit{-de\-formed magnetic field}:%
\begin{align}
B_{F} &  =\text{i}\hspace{0.01in}\partial^{C}\triangleright A^{D}%
\varepsilon_{DCF}=-q^{-4}\hspace{0.01in}(\hat{R}^{-1}){^{CD}}_{\hspace
{-0.01in}GH}\hspace{0.02in}A^{G}\triangleleft\partial^{H}\varepsilon
_{DCF}\nonumber\\
&  =\text{i}\hspace{0.01in}A^{G}\triangleleft\partial^{H}\varepsilon
_{HGF}=g_{FG}\hspace{0.02in}B^{G}.
\end{align}
Using the commutation relations in Eq.~(\ref{VerRelKinImp}), we can verify the
following identity by a direct calculation:%
\begin{equation}
\lbrack\Pi^{C}\Pi_{C},\Pi^{D}]=-2\text{i}\hspace{0.01in}m\hspace
{0.01in}f_{\text{Lor}}^{D}.\label{VerRelQuaKanImpKanImp}%
\end{equation}
Note that the $q$-\textit{de\-formed Lorentz force density}\textbf{ }reads as
follows:%
\begin{equation}
f_{\text{Lor}}^{D}=\frac{\text{i}e}{2\hspace{0.01in}m}\text{\hspace{0.01in}%
}g^{DG}\varepsilon_{ACG}\hspace{0.01in}(\Pi^{C}\hspace{0.01in}B^{A}%
-B^{C}\hspace{0.01in}\Pi^{A}).\label{LorKraDic}%
\end{equation}

\section{$q$-Analog of angular momentum\label{AppAngMom}}

We can express the components of $q$-de\-formed orbital angular momentum by
$q$-de\-formed partial derivatives and $q$-de\-formed position coordinates
\cite{Lorek:1995ph,Lorek:1997eh,Schmidt:2006}:%
\begin{equation}
L_{A}=q^{-2}\Lambda^{1/2}X^{C}\hat{\partial}^{D}\varepsilon_{DCA}%
=-\hspace{0.01in}q^{2}\Lambda^{1/2}\hspace{0.01in}\hat{\partial}^{C}%
\hspace{-0.01in}X^{D}\varepsilon_{DCA}.\label{BahDreDarOpe}%
\end{equation}
Once again, $\varepsilon_{DCA}$ denotes the $q$-de\-formed epsilon tensor [cf.
Eq.~(\ref{epseu3app}) in Chap.~\ref{AppLorFor}] and $\Lambda$ stands for a
scaling operator subject to the relations in Eq.~(\ref{VerSkaKooQuaEukDrei})
of Chap.~\ref{KapHofStr}. Together with the element%
\begin{equation}
W=\Lambda^{-1/2}\hspace{0.01in}(1+q^{3}\lambda\,g_{AB}\hspace{0.01in}%
X^{A}\partial^{B})=\Lambda^{1/2}\hspace{0.01in}(1+q^{-3}\lambda\,g_{AB}%
\hspace{0.01in}X^{A}\hat{\partial}^{B}),
\end{equation}
the three components of $q$-de\-formed angular momentum form the Hopf algebra
$\mathcal{U}_{q}(\operatorname*{su}_{2})$. In this respect, the following
relations hold \cite{Lorek:1995ph,Lorek:1997eh}:%
\begin{equation}
L^{A}L^{B}\varepsilon_{BAC}=Wg_{CD}L^{D}.
\end{equation}

{\normalsize
\bibliographystyle{unsrt}
\bibliography{book,habil}

\begin{thebibliography}{10}

\bibitem{Hagar:2014}
A. Hagar.
\newblock {\em Discrete or Continuous? The Quest for Fundamental Length in
  Modern Physics}.
\newblock Cambridge University Press, Cambridge/UK, 2014.

\bibitem{Heisenberg:1930}
W. Heisenberg.
\newblock Die {S}elbstenergie des {E}lektrons.
\newblock {\em Z. Phys.}, 64:4--63, 1930.

\bibitem{Heisenberg:1938}
W. Heisenberg.
\newblock Die {G}renzen der {A}nwendbarkeit der bisherigen {Q}uantentheorie.
\newblock {\em Z. Phys.}, 110:241--266, 1938.

\bibitem{Garay:1995}
{L. J.} Garay.
\newblock Quantum gravity and minimum length.
\newblock {\em Int. J. Phys. A}, 10:145--165, 1995.

\bibitem{Mead:1966zz}
C.~A. Mead.
\newblock Observable consequences of fundamental-length hypotheses.
\newblock {\em Phys. Rev.}, 143:990--1005, 1966.

\bibitem{CarowWatamura:1990nk}
U. Carow-Watamura, M. Schlieker, M.~Scholl, and S. Watamura.
\newblock Tensor representation of the quantum group {SL}${}_q${(2,C)} and
  quantum {M}inkowski space.
\newblock {\em Z. Phys. C}, 48:159--166, 1990.

\bibitem{Faddeev:1987ih}
L.~D. Faddeev, N.~{Yu}. Reshetikhin, and L.~A. Takhtajan.
\newblock Quantization of {L}ie groups and {L}ie algebras.
\newblock {\em Leningrad Math. J.}, 1:193--225, 1990.

\bibitem{Moyal:1949sk}
J.~E. Moyal.
\newblock Quantum mechanics as a statistical theory.
\newblock {\em Proc. Cambridge Philos. Soc.}, 45:99--124, 1949.

\bibitem{Carnovale:1999}
G. Carnovale.
\newblock On the braided {F}ourier transform in the $n$-dimensional quantum
  space.
\newblock {\em J. Math. Phys.}, 40:5972--5997, 1999.
\newblock arXiv:math/9810011.

\bibitem{Wachter:2007A}
H. Wachter.
\newblock Analysis on $q$-deformed quantum spaces.
\newblock {\em Int. J. Mod. Phys. A}, 22:95--164, 2007.
\newblock arXiv:math-ph/0604028.

\bibitem{Wachter:2020A}
H. Wachter.
\newblock Quantum dynamics on the three-dimensional $q$-deformed euclidean
  space.
\newblock 2020.
\newblock arXiv:math-ph/2004.05444.

\bibitem{Wachter:2020B}
H. Wachter.
\newblock Nonrelativistic one-particle problem on $q$-deformed euclidean space.
\newblock 2020.
\newblock arXiv:quant-ph/2010.08826.

\bibitem{Lorek:1997eh}
A. Lorek, W. Weich, and J. Wess.
\newblock Non-commutative {E}uclidean and {M}inkowski structures.
\newblock {\em Z. Phys. C}, 76:375--386, 1997.

\bibitem{Bayen:1977ha}
F.~Bayen, M.~Flato, C.~Fronsdal, A.~Lichnerowicz, and D.~Sternheimer.
\newblock Deformation theory and quantization. 1. {Deformations} of symplectic
  structures.
\newblock {\em Ann. Phys.}, 111:61--110, 1978.

\bibitem{1997q.alg.....9040K}
M.~{Kontsevich}.
\newblock Deformation quantization of {P}oisson manifolds, {I}.
\newblock arXiv:q-alg/9709040, 1997.

\bibitem{Madore:2000en}
J. Madore, S. Schraml, P. Schupp, and J. Wess.
\newblock Gauge theory on noncommutative spaces.
\newblock {\em Eur. Phys. J. C}, 16:161--167, 2000.
\newblock arXiv:hep-th/0001203.

\bibitem{Wachter:2002A}
H. Wachter and M. Wohlgenannt.
\newblock $*$-{P}roducts on quantum spaces.
\newblock {\em Eur. Phys. J. C}, 23:761--767, 2002.
\newblock arXiv:hep-th/0103120.

\bibitem{Jackson:1910yd}
F.~N. Jackson.
\newblock $q$-{D}ifference equations.
\newblock {\em Amer. J. Math.}, 32:305--314, 1910.

\bibitem{CarowWatamura:1990zp}
U. Carow-Watamura, M. Schlieker, and S. Watamura.
\newblock {SO}${}_{q}${($N$)}-covariant differential calculus on quantum space
  and quantum deformation of {S}chroe\-dinger equation.
\newblock {\em Z. Phys. C}, 49:439--446, 1991.

\bibitem{Wess:1990vh}
J. Wess and B. Zumino.
\newblock Covariant differential calculus on the quantum hyperplane.
\newblock {\em Nucl. Phys. Proc. Suppl. B}, 18:302--312, 1991.

\bibitem{Bauer:2003}
C. Bauer and H. Wachter.
\newblock Operator representations on quantum spaces.
\newblock {\em Eur. Phys. J. C}, 31:261--275, 2003.
\newblock arXiv:math-ph/0201023.

\bibitem{Wachter:2004A}
H. Wachter.
\newblock $q$-{I}ntegration on quantum spaces.
\newblock {\em Eur. Phys. J. C}, 32:281--297, 2004.
\newblock arXiv:hep-th/0206083.

\bibitem{Jackson:1908}
F.~N. Jackson.
\newblock On $q$-definite integrals.
\newblock {\em Quart. J. Pure and Appl. Math.}, 41:193--203, 1910.

\bibitem{Jambor:2004ph}
C. Jambor.
\newblock {\em Non-Commutative Analysis on Quantum Spaces}.
\newblock Dissertation, Fak. f. Phys., LMU M\"unchen, 2004.

\bibitem{Kempf:1994yd}
A. Kempf and S. Majid.
\newblock Algebraic $q$-integration and {F}ourier theory on quantum and braided
  spaces.
\newblock {\em J. Math. Phys.}, 35:6802--6837, 1994.
\newblock arXiv:hep-th/9402037.

\bibitem{Kulish:1983md}
P.~P. Kulish and N.~{\relax Yu}. Reshetikhin.
\newblock {Quantum linear problem for the Sine-Gordon equation and higher
  representation}.
\newblock {\em J. Sov. Math.}, 23:2435--2441, 1983.
\newblock [Zap. Nauchn. Semin.101,101(1981)].

\bibitem{Lorek:1993tq}
A. Lorek, W.~B. Schmidke, and J. Wess.
\newblock {SU}$_q(2)$-covariant {R}-matrices for reducible representations.
\newblock {\em Lett. Math. Phys.}, 31:279--288, 1994.

\bibitem{majid-1993-34}
S. Majid.
\newblock Braided momentum structure of the $q$-{P}oincar{\'e} group.
\newblock {\em J. Math. Phys.}, 34:2045--2058, 1993.
\newblock arXiv:hep-th/9210141.

\bibitem{Weixler:1993ph}
R.~O. Weixler.
\newblock {\em Inhomogene Quantengruppen}.
\newblock Dissertation, Fak. f. Phys., LMU M\"unchen, 1993.

\bibitem{Klimyk:1997eb}
A.~U. Klimyk and K. Schm\"udgen.
\newblock {\em Quantum groups and their representations}.
\newblock Springer, Berlin - Heidelberg - New York, 1997.

\bibitem{ogievetsky1992}
O.~Ogievetsky, W.~B. Schmidke, J.~Wess, and B.~Zumino.
\newblock $q$-deformed poincar\'e algebra.
\newblock {\em Comm. Math. Phys.}, 150(3):495--518, 1992.

\bibitem{Mikulovic:2006}
D. Mikulovic, A. Schmidt, and H. Wachter.
\newblock Grassmann variables on quantum spaces.
\newblock {\em Eur. Phys. J. C}, 45:529--544, 2006.
\newblock arXiv:hep-th/0407273.

\bibitem{Cerchiai:1999}
B.~L. Cerchiai, R. Hinterding, J. Madore, and J. Wess.
\newblock A calculus based on a $q$-deformed {H}eisenberg algebra.
\newblock {\em Eur. Phys. J. C}, 8:547--558, 1999.

\bibitem{Pauli:1980}
W. Pauli.
\newblock {\em General Principles of Quantum Mechanics}.
\newblock Springer, Berlin - Heidelberg - New York, 1980.

\bibitem{Lorek:1995ph}
A. Lorek.
\newblock {\em $q$-Deformierte Quantenmechanik und induzierte
  Wechselwirkungen}.
\newblock Dissertation, Fak. f. Phys., LMU M\"unchen, 1995.

\bibitem{2004JPhA...37.9175F}
G. Fiore.
\newblock Quantum group covariant (anti)symmetrizers, ${\varepsilon}$-tensors,
  vielbein, hodge map and {L}aplacian.
\newblock {\em J. Phys A}, 37:9175--9193, 2004.
\newblock arXiv:math/0405096.

\bibitem{majid-1995-36}
S. Majid.
\newblock $q$-{E}psilon tensor for quantum and braided spaces.
\newblock {\em J. Math. Phys.}, 36:1991--2007, 1995.
\newblock arXiv:hep-th/9406157.

\bibitem{Meyer:1994wi}
U. Meyer.
\newblock Wave equations on $q$-{M}inkowski space.
\newblock {\em Commun. Math. Phys.}, 174:457--476, 1995.
\newblock arXiv:hep-th/9404054.

\bibitem{Schmidt:2005ph}
A. Schmidt.
\newblock {\em $q$-Deformierte Superanalysis und Quanten-Lie-Algebren}.
\newblock Dissertation, Fak. f. Phys., LMU M\"unchen, 2005.

\bibitem{Schmidt:2006}
A. Schmidt and H. Wachter.
\newblock $q$-{D}eformed quantum {L}ie algebras.
\newblock {\em J. Geom. Phys.}, 56:2289--2325, 2006.

\end{thebibliography}
}
\end{document}